\documentclass[12pt]{article}
\usepackage{amsmath, amssymb, theorem, pgf, pgfplots, tikz, a4wide}

\usepackage{tikz}
	\usetikzlibrary{shapes,arrows}

\newtheorem{theorem}{Theorem}
\newtheorem{definition}[theorem]{Definition}\sl
\newtheorem{remark}[theorem]{Remark}\rm
\newtheorem{lemma}[theorem]{Lemma}
\newtheorem{proposition}[theorem]{Proposition}
\newtheorem{corollary}[theorem]{Corollary}
\newtheorem{ex}[theorem]{Example}

\newtheorem{conj}[theorem]{Conjecture}
\newtheorem{claim}[theorem]{Claim}
\newcommand{\bprf}{\bf Proof: \nopagebreak  \rm}
\newcommand{\eprf}{\nopagebreak\hspace*{1em} \hfill
$\Box$\noindent\vspace{2ex}\\\phantom{}} 

\usepackage{amsopn}

\DeclareMathOperator{\diag}{diag}
\DeclareMathOperator{\Real}{Re}
\DeclareMathOperator{\Imag}{Im}
\DeclareMathOperator{\spann}{span}
\DeclareMathOperator{\tr}{trace}

\begin{document}
  \title{Simultaneous hollowisation, joint numerical range,\\ and stabilization by noise}

\author{Tobias Damm\thanks{Department of
    Mathematics,  TU Kaiserslautern, Germany, \texttt{damm@mathematik.uni-kl.de}} \and Heike
  Fa\ss{}bender\thanks{Institut \emph{Computational Mathematics},
    TU Braunschweig, Germany, \texttt{h.fassbender@tu-bs.de}}}

\maketitle
\begin{abstract}
  We consider orthogonal transformations of arbitrary square matrices to a
  form where all diagonal entries are equal. In our main results we
  treat the simultaneous transformation of two matrices and
  the symplectic orthogonal transformation of one matrix. A relation
  to the joint real numerical range is worked out, efficient
  numerical algorithms are developped and applications to
  stabilization by rotation and by noise are presented.
\end{abstract}



\section{Introduction}
A square matrix whose diagonal entries are all zero, is sometimes called a
\emph{hollow matrix}, e.g.\ \cite{CharFarb13, FarbJohn15, KuraBapa16,
  NeveBast18}.  By a theorem of
Fillmore  \cite{Fill69}, which is closely related to older results of
Horn and Schur \cite{Schu23, Horn54}, every real square zero-trace matrix  is orthogonally similar to a hollow
matrix. Taken with a pinch of salt, the structure of a hollow matrix can be viewed as
the negative of the spectral normal form  (e.g.\ of a symmetric matrix), where the zeros are
placed outside the diagonal. While the spectral form reveals
an orthogonal basis of eigenvectors,  a hollow form reveals an
orthogonal basis of neutral vectors, i.e.\ vectors for which the
quadratic form associated to the matrix vanishes.\\
 This property turns out
to be relevant in asymptotic eigenvalue considerations. More concretely,
we use it to extend and give new proofs for results on stabilization
of linear systems by rotational forces or by noise. Since the
pioneering work \cite{ArnoCrau83} these phenomena have received
ongoing attention, with current interest e.g.\ in stochastic partial differential
equations or Hamiltonian systems,  \cite{SrinLali00, CaraRobi04, KolbConi19}. Our new contribution concerns simultaneous stabilization by
noise and features a new method of proof, which relies on
an orthogonal transformation of matrices to hollow form.\\ 
It is easy to see that -- in contrast to
the spectral transformation --  the
transformation to hollow form leaves a lot of freedom to require further properties.
 In the present note, we first show that it is
possible to transform two zero-trace matrices simultaneously to an almost hollow
form, as will be specified in Section \ref{sec:Hollow}.  In a non-constructive
manner,  the proof can be based on
Brickman's theorem \cite{Bric61} that the real joint numerical range of two real
matrices is convex.  But to make the transformation computable, we
provide a different proof, which is fully constructive. As a side result,
this also leads to a new derivation of Brickman's theorem.
Moreover, the simultaneous transformation result allows to prove a stronger version of
Fillmore's theorem, namely that every real square zero-trace matrix  is
orthogonal-symplectically similar to a hollow matrix. \\
We mainly treat the real case, because  it is slightly more
involved than its complex counterpart.   Complex versions of our results can be obtained easily and
are stated in Section \ref{sec:complex-case}. It turns out that any
pair of Hermitian zero-trace matrices is unitarily similar to a hollow
pair (not just an almost hollow pair as in the real case). 
In \cite{NeveBast18} the term \emph{simultaneous unitary
  hollowisation} is used for such a tranformation, and it is put in the context of a quantum
separability problem. The authors show that a certain quantum state is
separable, if and only if an associated set of matrices  is simultaneously unitarily
  hollowisable. This is a non-trivial restriction, if there are more
  than two matrices. 
For an arbitrary triple of Hermitian matrices, however, we can show
that it is  unitarily similar to an almost hollow
triple, i.e.\ \emph{almost hollowisable}, so to speak. We see this as a first
step towards criteria for larger sets of matrices to be simultaneously
unitarily (almost) hollowisable. Thus, to the best of our knowledge, the current note is the first to treat hollowisation
  problems from the matrix theoretic side.  \\
All our results are constructive and can be implemented in a
straight-forward way.  Computational aspects of the real  transformations are discussed in Section
\ref{sec:comp_aspects}. The orthogonal symplectic transformation of
$4\times 4$-matrices requires detailed explicit calculations which
have been shifted to the appendix.
We show analytically that for $n\times n$-matrices the
computational cost of our hollowising transformations is $O(n^2)$ and report on numerical
experiments. \\
In  Section \ref{sec:appl_stab}, we present the
applications of our results in stabilization theory. We show that a
number of  linear dissipative systems can be stabilized
simultaneously by the same stochastic noise process, provided the
coefficient matrices can be made almost hollow simultaneously by an
orthogonal transformation.  The results are illustrated by numerical examples.

\section{Hollow matrices and orthogonal transformations}
\label{sec:Hollow}
We first review some known facts on hollow matrices and then present our
main results.
\begin{definition}
  Let $A=(a_{ij})\in \mathbb{R}^{n\times n}$.
  \begin{enumerate}
  \item[(i)] We call $A$ \emph{hollow}, if $a_{ii}=0$ for all
    $i=1,\ldots,n$. \item[(ii)] We call $A$ \emph{almost hollow}, if $a_{ii}=0$
    for $i=1,\ldots,n-2$ and $a_{n-1,n-1}=-a_{nn}$. \item[(iii)] If $\tr A=0$, then
    $A$ is called a \emph{zero trace} matrix.
  \end{enumerate}

\end{definition}
Obviously, every hollow matrix is also almost hollow, and every almost
hollow matrix is zero trace. Vice versa, $\tr A=0$ implies that $A$ is orthogonally
similar to a hollow matrix. This result has been proven in Fillmore (1969)
\cite{Fill69}. We add a proof, because similar arguments will be used in
the later discussion.
\begin{lemma}\label{lemma:fillmore} Let $A\in\mathbb{R}^{n\times n}$ with $\tr A=0$.
  \begin{itemize}
  \item[(a)] There exists a vector $v\in\mathbb{R}^n$ with $v\neq 0$,
    such that $v^TAv=0$.
\item[(b)] There exists an orthogonal matrix $V\in\mathbb{R}^{n\times
    n}$, such that $V^TAV$ is hollow. 
  \end{itemize}
\end{lemma}
\bprf
(a) If $a_{11}=0$, then we can choose $v=e_1$. Otherwise let (after possibly
dividing $A$ by $a_{11}$) w.l.o.g.\
$a_{11}=1$. Since $\tr A=0$, there exists $j\in\{2,\ldots,n\}$ with
$a_{jj}<0$. For $v=xe_1+e_j$ with $x\in\mathbb{R}$, we have 
\begin{align*}
  v^TAv&=x^2+(a_{1j}+a_{j1})x+a_{jj}
\end{align*}
which has two real zeros.
Hence (a) follows.\\
(b)  Extend $v_1=v/\|v\|$ with $v$ from (a), to an orthonormal matrix
$V_1=[v_1,\ldots,v_n]$.  Then $V^TAV=\left[
  \begin{array}{cc}
    0&\star\\\star&A_1
  \end{array}
\right]$ with $A_1\in\mathbb{R}^{(n-1)\times (n-1)}$ and $\tr A_1=\tr A=0$.
Therefore we can proceed with $A_1$ as with $A$.
\eprf
\begin{corollary}\label{cor:fillmore}
 For $A\in\mathbb{R}^{n\times n}$, there exists an orthogonal matrix $V\in\mathbb{R}^{n\times
    n}$, such that all diagonal entries of $V^TAV$ are equal.
\end{corollary}
\bprf
We set $A_0=A-\frac{\tr A}n I$. By Lemma
\ref{lemma:fillmore} there exists an orthogonal matrix $V$ such that
$V^TA_0V$ is hollow. Then $V^TAV=V^TA_0V+\frac{\tr A}n I$.  
\eprf
\begin{remark}\label{rem:hollow}
  \begin{enumerate}
\item[(a)] A transformation matrix $V$ making $V^TAV$ hollow as in Lemma \ref{lemma:fillmore}
  will sometimes be called an \emph{(orthogonal) hollowiser (for $A$)}.
  \item[(b)] As is evident from the construction, the hollowiser
    $V$ is not unique. 
In the following
    we will exploit this freedom to transform two matrices
    simultaneously or to replace $V$ by an orthogonal symplectic
    matrix. 
\item[(c)] Since $V^TAV$ is hollow, if and only if $V^T(A+A^T)V$ is
  hollow, there is no restriction in considering only symmetric
  matrices.
\item[(d)] We are mainly interested in the real case, but it is
  immediate to transfer our results to the complex 
  case, where $A\in\mathbb{C}^{n\times n}$ and $V$ is unitary.  This
  is sketched in subsection \ref{sec:complex-case}.
  \end{enumerate}
\end{remark}

\subsection{Simultaneous transformation of two matrices}
\label{sec:simult-transf-two}
Simultaneous transformation of several matrices to a certain form (e.g.\ spectral
form) usually requires quite restrictive assumptions.  Therefore it is
remarkable that an arbitrary pair of zero trace matrices can simultaneously be
transformed to an almost hollow pair. The precise statement is given
in the following result.

\begin{proposition}\label{prop:fillmore_simultan}
  Consider $A,B\in\mathbb{R}^{n\times n}$ with
  $\tr A=\tr B=0$.
  \begin{itemize}
  \item[(a)] If $n\ge 3$,  there exists a nonzero vector $v\in\mathbb{R}^n$, such that $v^TAv=v^TBv=0$.
 \item[(b)] There exists an orthogonal matrix $V\in\mathbb{R}^{n\times
     n}$ such that $V^TAV$ is hollow and $V^TBV$ is almost hollow.
  \end{itemize}
\end{proposition}
\bprf
(b): We first note that (b) follows easily from (a). If (a) holds, then the orthogonal transformation $V$ is obtained by applying (a)
repeatedly as in the proof of Lemma \ref{lemma:fillmore}(b) until
the remaining submatrix is smaller than $3\times 3$ (where (a) is applied only for $A$). \\
For (a) we provide two different proofs. The first is quite short, but
not constructive. It exploits Brickman's theorem \cite{Bric61} on the
convexity of the joint real numerical
range of two matrices, see Theorem~\ref{lemma:realJNR_convex} below. The second is constructive, but considerably
longer. It is the basis for our algorithmic approach.

\emph{short proof} of (a):
By Lemma \ref{lemma:fillmore}, we can assume w.l.o.g.\ that $A$ is
hollow. If $b_{jj}=0$ for some $j$, then we
can choose $v=e_j$. Otherwise, since $\tr B=0$, not all the signs
of the $b_{jj}$ are equal. For simplicity of notation assume that
$b_{11}>0$ and  $b_{22}<0$. The points
$(e_1^TAe_1,e_1^TBe_1)=(0,b_{11})$ and $(e_2^TAe_2,e_2^TBe_2)
=(0,b_{22})$ lie in the joint real numerical range of $A$ and $B$,
defined as
 \begin{align*}
   W(A,B)&= \{(x^TAx,x^TBx)\;\big|\; x\in\mathbb{R}^n, \|x\|=1\}\subset\mathbb{R}^2\;.
  \end{align*}
According to Theorem~\ref{lemma:realJNR_convex} 
the set $W(A,B)$ is convex for $n\ge3$. Hence it also
contains $(0,0)=(v^TAv,v^TBv)$ for some unit vector $v\in\mathbb{R}^n$.

\emph{constructive proof} of (a): By Remark \ref{rem:hollow}, we can assume that $A$ and $B$ are
  symmetric, and by Lemma \ref{lemma:fillmore}, we can assume w.l.o.g.\ that $A$ is
hollow. If $b_{jj}=0$ for some $j$, then we
can choose $v=e_j$. For the remaining discussion let $b_{jj}\neq 0$ for all $j$.
Since $\tr B=0$, not all the signs of the $b_{jj}$ are equal. After possible
permutation and division of $B$ by one of the diagonal entries, we
can assume that the left upper $3\times 3$ blocks of $A$ and $B$ are
\begin{align}\label{eq:A3B3}
  A_3&=\frac12\left[\begin{array}{ccc}
    0&a&b\\a&0&c\\b&c&0
  \end{array}
\right]\;,\quad B_3=\frac12\left[\begin{array}{ccc}
    2d_-&\alpha&\beta\\\alpha&2d_+&\gamma\\\beta&\gamma&2
  \end{array}
\right]\,\quad \text{ with } d_-<0, d_+>0\;. 
\end{align}
If possible, we try to find $v_3=\left[
  \begin{smallmatrix}
    1\\x\\y
  \end{smallmatrix}
\right]$ with $x,y\in\mathbb{R}$, such that
$v_3^TA_3v_3=v_3^TB_3v_3=0$. This leads to the conditions
\begin{align}\label{eq:vAv}
  0&=v_3^TA_3v_3=ax+by+cxy=ax+(b+cx)y\\
0&=v_3^TB_3v_3=d_-+\alpha x+\beta y+\gamma xy+d_+x^2+y^2\;.
\label{eq:vBv}
\end{align}
We distinguish a number of cases.\\
{\bf\boldmath Case $a=0$ or $b=0$:}
If $a=0$, then \eqref{eq:vAv} holds with $y=0$ and \eqref{eq:vBv}
reduces to $0=d_-+\alpha x+d_+x^2$, which has a real solution $x$,
because $d_-<0$, $d_+>0$. Analogously, if $b=0$ then \eqref{eq:vAv}
holds with $x=0$ and \eqref{eq:vBv} again has a real solution. \\
{\bf\boldmath Case $a\neq0$, $b\neq0$, and $c\neq 0$:}
From now on let $a\neq 0$ and $b\neq 0$.
If equation \eqref{eq:vAv}  holds with $b+cx=0$, then also $ax=0$,
i.e.\ $a=0$ or $x=0$, where the latter implies $b=0$ and thus both
cases contradict our assumption. Therefore we can exclude the case
$b+cx=0$ and solve for $y=-\frac{ax}{b+cx}$. Inserting this in \eqref{eq:vBv} yields
\begin{align*}
  0&=d_-+\alpha x- \frac{\beta ax}{b+cx}-\frac{\gamma ax^2}{b+cx}+d_+x^2+\frac{a^2x^2}{(b+cx)^2}\;.
\end{align*}
If we multiply the equation with $(b+cx)^2$ and consider only the
coefficients at $x^0$ and $x^4$, we have
\begin{align}\label{eq:quartic}
  0&=d_-b^2+\ldots+d_+c^2x^4\;.
\end{align}
If $c\neq 0$, then $d_+c^2>0$ and $d_-b^2<0$ imply the existence of a
real root $x$. \\
{\bf\boldmath Case $a\neq0$, $b\neq0$, and $c= 0$:}
The final case to be considered is $c=0$. Now \eqref{eq:vAv} gives
$y=-\frac{a}bx$, which inserted in \eqref{eq:vBv} leads to
\begin{align*}
0&=  d_-+\alpha x-\beta \frac{a}bx+\left(-\gamma \frac{a}b+d_++\frac{a^2}{b^2}\right)x^2\;.
\end{align*}
Because $d_-<0$, the existence of a real root $x$ is guaranteed, if
 $-\gamma \frac{a}b+d_++\frac{a^2}{b^2}>0$. 

On the other hand note, that for any $\tilde v_3=\left[
  \begin{smallmatrix}
    0\\x\\y
  \end{smallmatrix}
\right]$, we have $\tilde v_3^T A_3\tilde v_3=0$. If moreover the
submatrix $\left[
  \begin{smallmatrix}
    2d_+&\gamma\\\gamma&2
  \end{smallmatrix}
\right]$ is not positive definite, i.e.\
$\gamma^2\ge 4d_+$, 
then there exists a
nonzero $\tilde v_3$, satisfying $\tilde v_3^T B_3\tilde v_3=0$. 

To conclude the proof, it suffices to note that the inequalities
$-\gamma \frac{a}b+d_++\frac{a^2}{b^2}\le 0$
and $\gamma^2< 4d_+$  contradict each other via
\begin{align*}
  0&\ge d_+-\gamma
     \frac{a}b+\frac{a^2}{b^2} > \frac{\gamma^2}4-\gamma
     \frac{a}b+\frac{a^2}{b^2}
     =\left(\frac\gamma2-\frac{a}b\right)^2\ge 0\;.
\end{align*}
The desired vector $v$ is now given by $v=\left[
  \begin{smallmatrix}
    v_3\\0
  \end{smallmatrix}
\right]$, or $v=\left[
  \begin{smallmatrix}
    \tilde v_3\\0
  \end{smallmatrix}
\right]$, respectively.
\eprf
\begin{remark}\rm
  The assumption in Proposition \ref{prop:fillmore_simultan}(a) that $n\ge
  3$ is essential. As the standard example (e.g.\ \cite{Bric61}), consider the two symmetric matrices $A=\left[
    \begin{array}{cc}
      1&0\\0&-1
    \end{array}
\right]$ and $B=\left[
    \begin{array}{cc}
      0&1\\1&0
    \end{array}
\right]$ with $\tr A=\tr B=0$. For $v= \left[
  \begin{smallmatrix}
    x\\y
  \end{smallmatrix}
\right]$, we have $v^TAv=x^2-y^2$ and $v^TBv=2xy$. If both forms are
zero, then necessarily $x=y=0$. Therefore, in general, a pair of symmetric
matrices with zero trace is not simultaneously orthogonally similar
to a pair of hollow matrices.
\end{remark}

\subsection{A constructive proof of Brickman's theorem}
\label{sec:constr-proof-brickm}

The following theorem was used in the short proof of Proposition \ref{prop:fillmore_simultan}(a).  It was derived in
\cite{Bric61, Bind85} by topological methods. More elementary
approaches using only connectivity properties of quadrics in
$\mathbb{R}^3$ were given in \cite{Yaku71, Pepi04, Mart05} and
surveyed e.g.\ in \cite{Poly98, PoliTerl07}.  
Below, we provide yet another
derivation, which exploits the $3\times 3$ case discussed in the
constructive proof.  While our approach might not be as elegant as some of the previous proofs, it easily lends itself for
computational purposes.
\begin{theorem}[Brickman \cite{Bric61}]\label{lemma:realJNR_convex}
  Let  $A,B\in\mathbb{R}^{n\times n}$ with $n\ge 3$. Then the set
  \begin{align*}
   W(A,B)&= \{(x^TAx,x^TBx)\;\big|\; x\in\mathbb{R}^n, \|x\|=1\}
  \end{align*}
is convex.
\end{theorem}
\bprf
Consider two linearly independent unit vectors $u,v\in\mathbb{R}^n$ and
set
\begin{align*}
a=(a_1,a_2)=(u^TAu,u^TBu), \;b=(b_1,b_2)=(v^TAv,v^TBv)\;.
\end{align*}
For $0<t<1$ let
$c= (c_1,c_2)=(1-t)a+tb$. We have to show that $c\in W(A,B)$, 
i.e.\
there exists a unit vector $x\in\mathbb{R}^n$, satisfying
$(x^TAx,x^TBx)=c$.\\
If $u^TAu=v^TAv$, then \emph{either}
$[u,v]^TA[u,v]=c_1I_2$, 
and we can choose $x\in\spann\{u,v\}$ ---
 \emph{or} $[u,v]^TA[u,v]-c_1I_2$ is indefinite, in which case there exist
 $z_\pm\in\spann\{u,v\}$ with $\|z_\pm\|=1$ such that $z_+^TAz_+>c_1$,
$z_-^TAz_-<c_1$. If  $u^TAu\neq v^TAv$, then we can trivially choose
$z_\pm\in\{u,v\}$ with the same properties. From now on, we assume
such vectors $z_\pm$ to be given.\\
Since $n\ge 3$ there exists another unit vector $y\in\mathbb{R}^n$
orthogonal to $z_\pm$. Depending on whether $y^TAy\ge c_1$ or
$y^TAy\le c_1$, we can choose a linear combination $w=\alpha y
+\beta z_-$ or $w=\alpha y
+\beta z_+$, $\alpha\neq 0$,  
such that $w^TAw=c_1$ and
$\|w\|=1$. 
With
the nonsingular matrix $U=[\sqrt{1-t}\,u,\sqrt{t}\,v,w]$, we define
\begin{align*}
  \tilde A&=U^T(A-c_1I)U=\left[
            \begin{array}{ccc}
             (1-t) (a_1-c_1)&\star&\star\\\star&t(b_1-c_1)&\star\\\star&\star&0
            \end{array}\right]\\
\tilde B&=U^T(B-c_2I)U=\left[
            \begin{array}{ccc}
             (1-t) (a_2-c_2)&\star&\star\\\star&t(b_2-c_2)&\star\\\star&\star&w^TBw-c_2
            \end{array}
\right]\;.
\end{align*}
By construction, $0=\tilde a_{11}+\tilde a_{22}=\tilde b_{11}+\tilde
b_{22}$. Hence, by Lemma \ref{lemma:fillmore}, there exists an orthogonal matrix
$Q_1\in\mathbb{R}^{2\times 2}$, such that for $Q=\left[
  \begin{array}{cc}
    Q_1&0\\0&1
  \end{array}
\right]$ we have
\begin{align*}
 Q^T \tilde A  Q&=\left[
            \begin{array}{ccc}
              0&\star&\star\\\star&0&\star\\\star&\star&0
            \end{array}\right],\quad
 Q^T \tilde B  Q=\left[
            \begin{array}{ccc}
              d_1&\star&\star\\\star&d_2&\star\\\star&\star&w^TBw-c_2
            \end{array}\right]\text{ where } d_1=-d_2\;.
\end{align*}
If 
$z^T Q^T \tilde A  Q z=z^T Q^T \tilde B  Q z=0$ for some vector
$z\in\mathbb{R}^3$, then $x=\frac{UQz}{\|UQz\|}\in\mathbb{R}^n$
yields
\begin{align*}
x^TAx&
       =\frac{z^TQ^T(\tilde A+c_1U^TU)Qz}{\|UQz\|^2}=c_1\;\text{ and }\; 
x^TBx
       =\frac{z^TQ^T(\tilde B+c_2U^TU)Qz}{\|UQz\|^2}=c_2,
\end{align*}
as desired.
Such a vector $z$ can be found as in the
constructive proof of Proposition~\ref{prop:fillmore_simultan}(a). If
$d_1=0$ or
$w^TBw=c_2$, then $x=e_1$ or $x=e_3$ is
suitable. Otherwise, after renormalization, the pair $(Q^T\tilde
AQ,Q^T\tilde BQ)$ has the same structure as $(A_3,B_3)$ in \eqref{eq:A3B3}.
This completes the proof.
\eprf


\subsection{Symplectic transformation of a matrix}
\label{sec:orth-sympl-transf}

Symplectic transformations play an important role in Hamiltonian
systems, e.g.\ \cite{MeyeHall09}. We briefly recapitulate some
elementary facts.
A real \emph{Hamiltonian matrix} has the form
\begin{align*}
  H&=\left[
     \begin{array}{cc}
       A&P\\Q&-A^T
     \end{array}
\right]\in\mathbb{R}^{2n\times 2n}\;,
\end{align*}
where $A \in\mathbb{R}^{n\times n}$ is arbitrary, while $P,Q
\in\mathbb{R}^{n\times n}$ are symmetric. If $J=\left[
     \begin{array}{cc}
       0&I\\-I&0
     \end{array}
\right]$, then all real Hamiltonian matrices are characterized by the
property that $JH$ is symmetric.  A real matrix $U \in\mathbb{R}^{2n\times 2n}$ is called
\emph{symplectic} if $U^TJU=J$.  If $U$ is symplectic, then the transformation $H\mapsto U^{-1}HU$ preserves the
Hamiltonian structure. Amongst other things, symplectic orthogonal
transformations are relevant for the Hamiltonian eigenvalue problem, e.g.\
\cite{PaigLoan81, Loan84, Fass00}. There is a rich theory on normal
forms of Hamiltonian matrices under orthogonal symplectic
transformations (e.g.\ \cite{Byer86, LinMehr99}). 
It is, however, a surprising improvement of Lemma \ref{lemma:fillmore} that an arbitrary zero trace matrix can
be  hollowised by a symplectic  orthogonal transformation. 
Before we state the main result of this section, we provide some
examples of symplectic orthogonal matrices, which will be relevant in
the proof and the computations.
\begin{ex}\rm
It is well-known and straight-forward to verify that an orthogonal matrix $U\in\mathbb{R}^{2n\times 2n}$
  is symplectic, if and only if it has the form
  \begin{align*}
U=\left[
    \begin{array}{cc}
      U_1&U_2\\-U_2&U_1
    \end{array}
  \right],
 \text{ where $U_1,U_2\in\mathbb{R}^{n\times n}$.}
  \end{align*}
This allows to construct elementary symplectic orthogonal
matrices (see e.g.\ \cite{MackMack03}).
\begin{enumerate}
\item If
  $V\in\mathbb{R}^{n\times n}$ is orthogonal, then $U=\left[
    \begin{smallmatrix}
      V&0\\0&V
    \end{smallmatrix}
  \right]$
  is symplectic orthogonal. 
\item If  $c^2+s^2=1$ then we define the Givens-type symplectic orthogonal matrices
  \begin{align}\label{eq:G2nics}
    G_k(c,s)&=\left[
        \begin{array}{ccc|ccc}
          I_{k-1}&&&&&\\&c&&&s&\\&&I_{n-k}&&&\\\hline&&&I_{k-1}&&\\&-s&&&c&\\&&&&&I_{n-k}
        \end{array}
\right],\; k\in\{1,\ldots,n\}\\
\mathcal{G}(c,s)&= 
                                      \left[
        \begin{array}{ccc|ccc}
          I_{n-2}&&&&&\\&c&s&&&\\&-s&c&&&\\\hline&&&I_{n-2}&&\\&&&&c&s\\&&&&-s&c
        \end{array}
\right]\;. \label{eq:calGncs}
\end{align}
\item For $p_0^2+p_1^2+p_2^2+p_3^2=1$ we have the symplectic
  orthogonal $4\times 4$-matrix
\begin{align}
S &= \begin{bmatrix}
p_0 & -p_1 & -p_2& -p_3\\
p_1 & p_0 & -p_3 & p_2\\
p_2 & p_3 & p_0 & -p_1\\
p_3 & -p_2 & p_1 & p_0
\end{bmatrix}\;.
\label{eq:4x4sympl}
 \end{align}

\end{enumerate}

\end{ex}\rm
\begin{theorem}\label{thm:SymplOrth}
   Consider a  matrix $A\in\mathbb{R}^{2n\times 2n}$ with 
   $n\ge1$. Then there exists a symplectic orthogonal
   matrix $U$, such that $U^TAU$ has constant diagonal.
\end{theorem}

\bprf
W.l.o.g.\ we can assume that $A$ is symmetric with $\tr A=0$. 
The transformation $U$ is constructed in several steps, where we make
use of the orthogonal symplectic transformations above.

{\bf 1st step:}
Let $d_1,\ldots,d_{2n}$ denote the diagonal entries of $A$.
Applying $G_k(c,s)$ from \eqref{eq:G2nics} for the transformation $A^+=G_k(c,s)^TA G_k(c,s)$
we can achieve that $d_{k}^+=d_{k+n}^+$.
After $n$ such transformations we have
\begin{align}
A^+&=\left[
  \begin{array}{cc}A_{1}^+&\star\\\star&A_2^+ \end{array}
\right]=
\left[  \begin{array}{cc}
    \begin{array}{ccc}
      d_{1}^+&&\star\\&\ddots&\\\star&&d_{n}^+
    \end{array}
    &\star\\\star&\begin{array}{ccc}
      d_{1}^+& &\star\\&\ddots&\\\star&&d_{n}^+
    \end{array}
  \end{array}
\right]\;.\label{eq:Aplus}
\end{align}
In particular $\tr A_1^+=\tr A_2^+=0$. 

{\bf 2nd step:} By Proposition
\ref{prop:fillmore_simultan}, there exists an orthogonal matrix
$V\in\mathbb{R}^{n\times n}$, such that $V^TA_1^+V$ is hollow and
$V^TA_2^+V$ is almost hollow.
Thus, for the symplectic orthogonal matrix $U=\left[
  \begin{array}{cc}
    V&0\\0&V
  \end{array}
\right]$, we have (with $d_1=0$)
\begin{align*}
U^TA^+U&=
  \left[
  \begin{array}{c|c}V^TA_{1}^+V&\star\\\hline\star&V^TA_2^+V \end{array}
\right]=
\left[  \begin{array}{c|c}
    \begin{smallmatrix}
      0&&\star\\[-2mm]&\ddots&\\\star&&\left[
                                  \begin{smallmatrix}
                                    d_1&a\\a&-d_1
                                  \end{smallmatrix}
\right]
    \end{smallmatrix}
    &\star\\\hline\star&\begin{smallmatrix}
      0& &\star\\[-2mm]&\ddots&\\\star&&\left[
                                  \begin{smallmatrix}
                                    d_2&b\\b&-d_2
                                  \end{smallmatrix}
\right]
    \end{smallmatrix}
  \end{array}
\right]\;.
\end{align*}

{\bf 3rd step:}
In the following we can restrict our attention to the submatrix of
$U^TA^+U$ formed by the rows and columns with indices $n-1,n,2n-1,2n$.
Therefore, we now work with symplectic orthogonal matrices
$G_k(c,s)$ from \eqref{eq:G2nics}, where $k\in\{n-1,n\}$ or
$\mathcal{G}(c,s)$ from \eqref{eq:calGncs}.  
Then it suffices to transform a $4\times 4$ symmetric matrix $A_4= \left[
  \begin{array}{cccc}
    d_1&a&\star&\star\\a&-d_1&\star&\star\\\star&\star&d_2&b\\\star&\star&b&-d_2
  \end{array}
\right]$ 
with the symplectic Givens rotations

\begin{align*} 
 G_{12}&=
\left[\begin{array}{cccc}
  c&s&0&0\\-s&c&0&0\\0&0&c&s\\0&0&-s&c
\end{array}\right],
\; G_{13}=
\left[\begin{array}{cccc}
  c&0&s&0\\0&1&0&0\\-s&0&c&0\\0&0&0&1
\end{array}\right],\;
G_{24}=\left[\begin{array}{cccc}
1&0&0&0\\0&  c&0&s\\0&0&1&0\\0&-s&0&c
\end{array}\right]\;.
\end{align*}
In an iterative approach, we show that for each such matrix $A_4$ with $|d_1|+|d_2|\neq0$ there
exists a product $G$ of matrices from this list so that 
\begin{align*}
  G^TA_4G=\left[
  \begin{array}{cccc}
    d_1^+&a^+&\star&\star\\a^+&-d_1^+&\star&\star\\\star&\star&d_2^+&b^+\\\star&\star&b^+&-d_2^+
  \end{array}
\right]\quad\text{ with } \quad |d_1^+|+|d_2^+| < |d_1|+|d_2|\;.
\end{align*}
We distinguish between different cases.

If $d_1\neq d_2$, then we can apply transformations with suitable $G_{13}$
and $G_{24}$ so that $d_1^+=d_2^+=(d_1+d_2)/2$. In particular $
 |d_1^+|+|d_2^+| \le |d_1|+|d_2|$.

If $d_1=d_2=:d$, let us assume w.l.o.g.\ that $|a|\ge|b|$. Moreover
assume that $d>0$ and $a>0$. Other combinations can be treated
analogously to the following considerations.
Setting
$A_4^+=G_{12}^TA_4G_{12}$, we have
\begin{align*}
d_1^+& =d_1^+(c,s)
       =   d(c^2-s^2) - 2acs \;,\quad
d_2^+ =d_2^+(c,s)
       =   d(c^2-s^2) - 2bcs \;.
\end{align*}
If $c=\cos(t)$,  $s=\sin(t)$, then $d_1^+$ is positive for $t=0$, 
negative for $t=\pi/4$ and strictly decreasing in $t$ on the interval $[0,\pi/4]$.
A direct calculation shows that $d_1^+=0$ for 
\begin{align}\label{eq:defcs}
  c&=\left(\frac12+\frac{a}2(d^2+a^2)^{-1/2}\right)^{1/2}\;,\quad
s=\left(\frac12-\frac{a}2(d^2+a^2)^{-1/2}\right)^{1/2}\;.
\end{align}
Here $c=\cos(t_0)$, $s=\sin(t_0)>0$ with minimal $t_0\in ]0,\pi/4[$,
and therefore  $c^2>s^2$ and $c,s>0$. 
Hence, if $b\ge 0$, then $a\ge b$ implies
\begin{align*}
  d&>  d_2^+= d(c^2-s^2) - 2bcs \ge 0\;.
\end{align*}
In this case $|d_1^+|+|d_2^+|=|d_2^+| \le |d|=\frac12
(|d_1|+|d_2|)$ as desired.\\
The case $b<0$ is slightly more subtle. We first derive a lower bound
for $s$ in \eqref{eq:defcs}. To this end note that the norm $\|A_4\|_2=\Delta$ is
invariant under orthogonal transformations and $a\le\Delta$.
Hence, for a given $d>0$, we have  
$$s\ge \left(\frac12-\frac{\Delta}2(d^2+\Delta^2)^{-1/2}\right)^{1/2}=:\mu(d)>0\;.$$
Since $d_1^+,d_2^+\ge 0$, we have
\begin{align*}
|d_1^+|+|d_2^+|=  d_2^+(c,s)&=2d(c^2 -s^2) - (a+b)cs\le 2d(1-2s^2)\le 2d(1-\mu(d)^2)<2d\;.
\end{align*}

Altogether, given $A_4$ we set $G=G_{12}G_{13}G_{24}$, where the
transformation with $G_{13}G_{24}$ achieves $d_1=d_2$ and $G_{12}$
makes $|d_1|+|d_2|$ smaller. Applying these transformations repeatedly, we obtain a sequence
$[d_1^{(k)},d_2^{(k)}]$ of diagonal entries, whose norm $|d_1^{(k)}|+|d_2^{(k)}|$ is monotonically
decreasing, and in the limit necessarily $\mu(d)=0$, which implies that
$[d_1^{(k)},d_2^{(k)}]\stackrel{k\to\infty}\longrightarrow 0$.
\eprf

\begin{remark}
  The previous proof is constructive, but the iterative approach to
  the $4\times 4$ case in the 3rd step is numerically inefficient. In
  the appendix we provide a direct construction of the transformation,
  which exploits also transformations of the special type \eqref{eq:4x4sympl}. 
\end{remark}

\subsection{The complex Hermitian case}
\label{sec:complex-case}
 The joint numerical range has been studied in even more detail for the
 complex Hermitian case than for the real case. Some of our results
 simplify or become even stronger if we allow for complex  unitary 
 instead of real orthogonal transformations. In the current subsection
 we sketch briefly how the results can be transferred.
For completeness we start with the complex version of Lemma
\ref{lemma:fillmore}, whose immediate proof is omitted, see \cite{Fill69}.
 \begin{lemma}\label{lemma:fillmore_complex} Let
   $A\in\mathbb{C}^{n\times n}$ be Hermitian with $\tr A=0$.
  Then there exists a unitary matrix $V\in\mathbb{C}^{n\times
    n}$, such that $V^*AV$ is hollow. 
\end{lemma}
From our approach it is less obvious than in the real case 
that the statement of this lemma holds for non-Hermitian $A$, too (a
fact already proven in \cite{Fill69}).
Our proof of Lemma \ref{lemma:fillmore} requires realness of the
diagonal entries, and in contrast to Remark \ref{rem:hollow}(c), the
property of $V^*AV$ being hollow is not equivalent to $V^*(A+A^*)V$
being hollow (take e.g.\ $A=iI$).
We will obtain the non-Hermitian version of Lemma
\ref{lemma:fillmore_complex} as a consequence of Proposition
\ref{prop:fillmore_simultan_complex} below. For the other statements
in this subsection we are not able to drop the Hermitian assumption
(see also Remark \ref{rem:Counterexamples}).
A complex version of Brickman's theorem has been
proven in \cite{Bind85}.
\begin{theorem}\label{thm:complexJNR_convex}
  Consider Hermitian matrices $A,B,C\in\mathbb{C}^{n\times n}$. Depending on $n$, the
  following sets are convex:
  \begin{align*}
n\ge1:&   \quad W(A,B):= \{(x^*Ax,x^*Bx)\;\big|\; x\in\mathbb{C}^n, \|x\|=1\}\;,\\
 n\ge 3:&  \quad W(A,B,C):= \{(x^*Ax,x^*Bx,x^*Cx)\;\big|\; x\in\mathbb{C}^n, \|x\|=1\}\;.
  \end{align*}
\end{theorem}
Based on Theorem~\ref{thm:complexJNR_convex}, it is easy to derive complex versions of Proposition
\ref{prop:fillmore_simultan} and Theorem~\ref{thm:SymplOrth}. 
\begin{proposition}\label{prop:fillmore_simultan_complex}
Let $A,B,C\in\mathbb{C}^{n\times n}$ be zero-trace Hermitian matrices.
  \begin{itemize}
  \item[(a)] If $n\ge 3$,  there exists  $v\in\mathbb{C}^n\setminus\{0\}$, such that $v^*Av=v^*Bv=v^*Cv=0$.
 \item[(b)] There exists a unitary matrix $V\in\mathbb{C}^{n\times
     n}$ such that $V^*AV$ and $V^*BV$ are hollow, while $V^*CV$ is almost hollow.
  \end{itemize}
\end{proposition}
\bprf
  We first consider only $A$ and $B$. By Lemma
  \ref{lemma:fillmore_complex}, we can assume $A$ to be hollow.  Literally
  as in the short proof of Proposition \ref{prop:fillmore_simultan}(a)
  it follows then that $(0,0)$ lies in the convex hull of $W(A,B)$ and
  thus in $W(A,B)$ itself by Theorem~\ref{thm:complexJNR_convex}.
  Hence, as in the proof of
  Proposition \ref{prop:fillmore_simultan}(b), there exists a unitary
  matrix $V$, such that $V^*AV$ and $V^*BV$ are hollow. If $n<3$ this
  proves (b).\\
 If $n\ge 3$ we assume for simplicity that $A$ and $B$ are already
 hollow. If one of the diagonal entries of $C$ vanishes, say
 $c_{jj}=0$, then we can choose $v=e_j$. Otherwise, there exist
 $j,k\in\{1,\ldots,n\}$ such that $c_{jj}c_{kk}<0$. Since
 $(0,0,c_{jj}),(0,0,c_{kk})\in W(A,B,C)$, another application of Theorem~\ref{thm:complexJNR_convex} yields $0\in W(A,B,C)$ and thus
 (a). As before, (b) is a consequence of (a).
\eprf

\begin{corollary}\label{cor:fillmore_complex}
  Let   $A\in\mathbb{C}^{n\times n}$ with $\tr A=0$.
  Then there exists a unitary matrix $V\in\mathbb{C}^{n\times
    n}$, such that $V^*AV$ is hollow. 
\end{corollary}
\bprf
  The matrices $\Real A=\frac12(A+A^*)$ and $\Imag A=\frac1{2i}(A-A^*)$ are
  Hermitian with zero trace. By Proposition
  \ref{prop:fillmore_simultan_complex}(b), there exists a unitary $V$
  such that $V^*(\Real A) V$ and $V^*(\Imag A) V$ are hollow. Thus $V^*AV$ is
  hollow as well.
\eprf
\begin{corollary}\label{cor:SymplUnit}
   Consider a  Hermitian matrix $A\in\mathbb{C}^{2n\times 2n}$ with
   $\tr A=0$.
   \begin{itemize}
   \item[(a)] There exists a unitary matrix $U$, such that $U^*JU=J$
     and $U^*AU$ is hollow.
  \item[(b)] There exists a unitary matrix $U$, such that $U^TJU=J$
     and $U^*AU$ is hollow.
   \end{itemize}
\end{corollary}
In the terminology of \cite{MackMack03}, the unitary matrix $U$ is
called \emph{conjugate
symplectic} in (a) and \emph{complex symplectic} in (b).\\
\bprf
  We repeat the first two steps in the proof of Theorem
  \ref{thm:SymplOrth}. Since $A$ is Hermitian, the first step can be
  carried out with a real transformation. Therefore we can assume that
  $A$ has the form $A=A^+$ from \eqref{eq:Aplus}. By Proposition
  \ref{prop:fillmore_simultan_complex}, there exists a unitary
  $V\in\mathbb{C}^{n\times n}$ such that $V^*A_1^+V$ is hollow and (a)
  $V^*A_2^+V$ is hollow or (b) $V^*\bar A_2^+V$ is hollow. Then $U=\left[
    \begin{smallmatrix}
      V&0\\0&V
    \end{smallmatrix}
\right]$ fulfils (a) or $U=\left[
    \begin{smallmatrix}
      V&0\\0&\bar V
    \end{smallmatrix}
\right]$ fulfils (b), respectively. 
\eprf
\begin{remark}\label{rem:Counterexamples}
If some of the assumptions are dropped, we can produce counter\-examples
to the statements of 
  Theorem \ref{thm:complexJNR_convex} and Proposition
  \ref{prop:fillmore_simultan_complex}.
  \begin{enumerate}
  \item Let $n=2$. For Hermitian matrices $A,B,C\in\mathbb{C}^{2\times 2}$ the
    set $W(A,B,C)$ needs not be convex. In \cite{GutkJonc04} it was
    shown that $W(A,B,C)$ is the unit sphere in $\mathbb{R}^3$ for $$A=\left[
      \begin{array}{cc}
        1&0\\0&-1
      \end{array}
\right]\;,\quad B=\left[
      \begin{array}{cc}
        0&1\\1&0
      \end{array}
\right]\;,\quad C=\left[
      \begin{array}{cc}
        0&i\\-i&0
      \end{array}
\right]\;.$$
In particular $W(A,B,C)$ is not convex and $0\not\in W(A,B,C)$ for these matrices, implying that
there is no $v\neq0$ with $v^*Av=v^*Bv=v^*Cv=0$. 
\item For non-Hermitian zero-trace matrices $A,B\in\mathbb{C}^{n\times
    n}$ there might be no $v\neq0$ with $v^*Av=v^*Bv=0$, and
  (consequently) $W(A,B)$ may
  be non-convex. As an example for arbitrary $n\ge 2$ consider
  \begin{align*}
    A=\left[
    \begin{array}{cc|c}
      1&0&0\\0&-1-(n-2)i&0\\\hline 0&0&iI_{n-2}
    \end{array}
\right]\;,\quad B=\left[
    \begin{array}{cc|c}
      0&0&0\\1&0&0\\\hline0&0&0_{n-2}
    \end{array}
\right] \;.
  \end{align*}
Obviously $(e_1^*Ae_1,e_1^*Be_1)=(1,0)$ and
$(e^*Ae,e^*Be)=\left(-(n-1)^{-1/2},0\right)$ for $e=(n-1)^{-1/2}\sum_{j=2}^ne_j$
with $\|e\|=1$. Hence $(0,0)$ lies in the convex hull of $W(A,B)$. But
the ansatz $(v^*Av,v^*Bv)=(0,0)$ with 
$v=\left[\begin{smallmatrix}
    x\\y\\z
  \end{smallmatrix}\right]$, $x,y\in\mathbb{C}$, $z\in\mathbb{C}^{n-2}$ yields
\begin{align*}
  |x|^2-|y|^2+i\big(\|z\|^2-(n-2)|y|\big)=0\;\text{ and }\; \bar y x=0\;.
\end{align*}
By the second equation we have $x=0$ or $y=0$. Together with the real part of
the first equation this implies $x=y=0$. The imaginary part of the
first equation then yields also $z=0$, i.e.\ $v=0$.  
  \end{enumerate}
\end{remark}

\section{Computational aspects}
\label{sec:comp_aspects}
The orthogonal transformation of a single matrix $A$ with $\tr A=0$
to a hollow matrix is straightforward along the proof of Lemma
\ref{lemma:fillmore}. Note that each nonzero diagonal entry can be
eliminated by one Givens rotation. Hence, if there are $\nu$ nonzero
diagonal entries, then $\nu-1$ Givens rotations are required.

\subsection{Simultaneous transformation of two matrices}
The transformation of a pair $(A,B)$ of zero trace matrices follows
the constructive proof of Proposition \ref{prop:fillmore_simultan}. 
In the first step, $A$ is transformed to hollow form.

Given a pair $(A,B)$ of $k\times k$ matrices, $k\ge 3$, with $A$ hollow and $B$ zero-trace, we first
check, whether  $b_{11}=0$. If so, then the
dimension can be reduced immediately.

Else, let $i_2\neq i_3$ with
$b_{i_2,i_2}=\min\{b_{22},\ldots,b_{nn}\}$ and  $b_{i_3,i_3}=\max\{b_{22},\ldots,b_{nn}\}$. For the submatrices $A_3$ of $A$
and $B_3$ of $B$ corresponding to the rows and columns $1,i_2,i_3$ as in \eqref{eq:A3B3}, a common neutral vector $v_3\in\mathbb{R}^3$
is computed. Generically, this requires the solution of a quartic
equation as in \eqref{eq:quartic}. The vector $v_3$ can be extended to
an orthogonal $k\times k$ matrix $V$ which differs from a permutation
matrix only in a $3\times 3$ subblock. After the transformation
\begin{align}
(A,B)\leftarrow
(V^TAV,V^TBV)=\left(\left[
  \begin{array}{c|c}
    0&\star\\\hline \star &\tilde A
  \end{array}
\right], \left[\begin{array}{c|c}
    0&\star\\\hline\star&\tilde B
  \end{array}
\right]\right)\label{eq:trafoAB}
 \end{align}
we have $\tr\tilde A=\tr\tilde B=0$, where at most
two diagonal entries of $\tilde A$ are non-zero. Hence, by another Givens rotation we have reduced the problem from dimension $k$
to  $k-1$. Since each Givens rotation and each transformation
\eqref{eq:trafoAB}
 requires $O(k)$ elementary operations, the whole algorithm
has complexity $O(n^2)$ including the solution of at most $n-2$
quartic equations.\\
We carried out experiments on a 2016 MacBook Pro with a 3.3 GHz
Intel Core i7 processor and 16 GB Memory running OS X 10.14.6 using
MATLAB version R2019b. For $20$ random pairs of $n\times n$ matrices $A$, $B$ we
avaraged the computing times, see Table \ref{tab:comp_sim}. Although the theoretical complexity is
not manifest in  the outcome, we see that the algorithm
is quite fast also for large matrices.
\begin{table}[tbhp]
  \centering
  \begin{tabular}{c|rrrrrrr}
   size $n$&100&200&400&800&1600&3200&6400\\\hline
    time in sec&0.016&0.037&0.10&0.72&7.9&71&852
  \end{tabular}
\caption{Computing times for simultaneous orthogonal transformation to
hollow form}
\label{tab:comp_sim}
\end{table}
\subsection{Symplectic transformation of a matrix}
The symplectic orthogonal transformation of a single matrix follows
the three steps in the proof of Theorem~\ref{thm:SymplOrth}. In the
3rd step the direct construction in Appendix \ref{app.1} is used.
This also gives an algorithm of complexity $O(n^2)$. Numerical
experiments with MATLAB were carried out as in the previous
subsection.  Again, the theoretical complexity is not really expressed
by the computing times in Table \ref{tab:comp_sym} (or only roughly between $2n=200$ and $2n=800$), but
most likely this is due to other effects such as memory management for large $n$.
\begin{table}[tbhp]
  \centering
  \begin{tabular}{c|rrrrrrr}
   size $2n$&100&200&400&800&1600&3200&6400\\\hline
    time in sec&0.010 &   0.013&    0.038 &   0.17  &  1.2&   11.7 &  97
  \end{tabular}
\caption{Computing times for symplectic orthogonal transformation to
hollow form}
\label{tab:comp_sym}
\end{table}
\section{Applications to stabilization problems}
\label{sec:appl_stab}
In this section we present two related stabilization problems. Both
deal with unstable linear ordinary differential equations, whose coefficient
matrix has negative trace. Such systems have stable and unstable
modes, but the stable ones dominate. By a mixing of the modes the system
can be stabilized. This mixing can be  achieved e.g.\ by adding
rotational forces or stochastic terms. For both cases we extend known
results from the literature. The basic idea lies in an asymptotic
analysis based on the hollow forms constructed in the previous sections. 

\subsection{Hamiltonian stabilization by rotation}
\label{sec:hamilt-gyrosc-stab}
A linear autonomous system
$
  \dot x=Ax
$
is called asymptotically stable, if all solutions $x(t)$ converge to
$0$ for $t\to\infty$. It is well known, that this is equivalent to the spectrum of
$A$ being contained in the open left half plane,
$\sigma(A)\subset\mathbb{C}_-$. In this case, necessarily $\tr
A<0$. Vice versa, one can ask, whether for any matrix $A$ with $\tr
A<0$, there exists a zero trace matrix $M$ of a certain type, such that
$\sigma(A+M)\subset\mathbb{C}_-$. In \cite{CrauDamm07} it
has been shown, that such a matrix $M$ can always be chosen to be
skew-symmetric. Then we say that $M$ stabilizes $A$ or
by rotation, see e.g.\ \cite{BaxeHenn93}.
The following theorem extends this result.
\begin{theorem}\label{thm:HamSkew}
  Let $A\in\mathbb{R}^{2n\times 2n}$ with $\tr A<0$. Then there
  exists a skew-symmetric Hamiltonian matrix $M$, such that $\sigma(A+M) \subset\mathbb{C}_-$.
\end{theorem}
\bprf
By Theorem~\ref{thm:SymplOrth} there exists a symplectic orthogonal matrix $U$, such
that $U^TAU$ has all diagonal entries equal to $\alpha=\tfrac{\tr A}{2n}<0$.
Consider $M_0=\left[
  \begin{array}{cc}
   0&\Lambda\\-\Lambda&0
  \end{array}
\right]$ with
$\Lambda=\diag(\lambda_1,\ldots,\lambda_n)\in\mathbb{R}^{n\times n}$, where
$|\lambda_j|\neq|\lambda_k|$ for $j\neq k$. Then $M_0$ is Hamiltonian and
skew-symmetric, and  
has all simple eigenvalues
$\pm i\lambda_k$
with respective eigenvectors $e_k\pm ie_{k+n}$. \\
For $\varepsilon>0$ we perturb $M_0$ to $M_\varepsilon=M_0+\varepsilon U^TAU$.
By \cite[Theorem 3.1]{CrauDamm07} (see also \cite{StewSun90,
  HinrPrit05}) the eigenvalues of $M_\varepsilon$ have the expansion
\begin{align*}
\pm i\lambda_k+\varepsilon(e_k\pm ie_{k+n})^*U^TAU (e_k\pm ie_{k+n})+O(\varepsilon^2)
  =\pm i\lambda_k+\varepsilon\alpha+O(\varepsilon^2)\;.
\end{align*}
Hence 
\begin{align*}
  \sigma(A+\tfrac1\varepsilon UM_0U^T)=\{\alpha\pm\tfrac1\varepsilon i\lambda_k+O(\varepsilon)\;\big|\;k=1,\ldots,n\}\subset\mathbb{C}_- 
\end{align*}
for sufficiently small $\varepsilon$. The matrix $M=\frac1\varepsilon
UM_0U^T$ stabilizes $A$ by rotation. Since $U$ is symplectic
orthogonal, the matrix $M$ is skew-symmetric Hamiltonian.
\eprf
\begin{ex}\rm
  We illustrate Theorem~\ref{thm:HamSkew}  by $A=\diag(1,1,1,-4)$
  and $M_0$ as above with $\Lambda=\diag(1,2)$. The matrix $A$ is
  hollowised by the orthogonal symplectic matrix
  $U=\tfrac12\left[
    \begin{smallmatrix} \sqrt{2} & \sqrt{2} & 0 & 0\\ 1 & -1 & 1 &
      -1\\ 0 & 0 & \sqrt{2} & \sqrt{2}\\ -1 & 1 & 1 & -1
    \end{smallmatrix}\right]$. Then $\tilde M_0=UM_0U^T=\tfrac14\left[
    \begin{smallmatrix} 
0 & -\sqrt{2} & 6 & -\sqrt{2}\\ \sqrt{2} & 0 & -\sqrt{2} & 6\\ -6 &
\sqrt{2} & 0 & -\sqrt{2}\\ \sqrt{2} & -6 & \sqrt{2} & 0 
    \end{smallmatrix}\right]$ is skew-symmetric and Hamiltonian. The spectral abscissa
  $\alpha(\mu)=\max\Real\sigma\left(A+\mu \tilde M_0\right)$ for $\mu>0$ is depicted in Fig.\
  \ref{fig:specAbscHam}. It becomes negative for $\mu\approx3.7$.
  Hence for $\mu>3.7$ the system $\dot x=(A+\mu \tilde M_0)x$ is
  asymptotically stable. In \cite{CrauDamm07} a servo-mechanism was
  described, which chooses a suitable gain $\mu$ adaptively via the
  feedback equation
  \begin{align}\label{eq:adaptively}
    \dot x&=(A+\mu(t) \tilde M_0)\,x\;,\quad \dot\mu =\|x(t)\|\;.
  \end{align}
  This method also works in the current example (see the right plot in
  Fig.\  \ref{fig:specAbscHam}), where $\mu$ roughly
  converges to $e^{2.73}-1\approx14.37$.

  \begin{figure}[h]\centering
    \begin{minipage}{.48\linewidth}
%
%
\definecolor{mycolor1}{rgb}{0.00000,0.44700,0.74100}%
\begin{tikzpicture}

\begin{axis}[%
width=5.2cm,
height=3.5cm,
at={(0cm,0cm)},
scale only axis,
xmin=0,
xmax=20,
ymin=-0.4,
ymax=1,
axis background/.style={fill=white},
xmajorgrids,
ymajorgrids,
legend style={legend cell align=left, align=left, draw=white!15!black}
]
\addplot [color=mycolor1, line width=1.5pt]
  table[row sep=crcr]{%
0	1\\
0.202020202020202	0.996323449769375\\
0.404040404040404	0.98525224947444\\
0.606060606060606	0.96666121280191\\
0.808080808080808	0.94034043489865\\
1.01010101010101	0.905995615106584\\
1.21212121212121	0.863254439244406\\
1.41414141414141	0.811689349457176\\
1.61616161616162	0.750878815168261\\
1.81818181818182	0.680551914294432\\
2.02020202020202	0.600895940918772\\
2.22222222222222	0.513123003548879\\
2.42424242424242	0.420246728591035\\
2.62626262626263	0.327471011218999\\
2.82828282828283	0.241050765838178\\
3.03030303030303	0.165719708531341\\
3.23232323232323	0.103070046553813\\
3.43434343434343	0.0521795679426701\\
3.63636363636364	0.0110818477036396\\
3.83838383838384	-0.0222493511348125\\
4.04040404040404	-0.0495245582977444\\
4.24242424242424	-0.0720815753114523\\
4.44444444444444	-0.0909378251419724\\
4.64646464646465	-0.106862283420873\\
4.84848484848485	-0.120438430299843\\
5.05050505050505	-0.132112797190283\\
5.25252525252525	-0.142230704122053\\
5.45454545454545	-0.151062158462698\\
5.65656565656566	-0.158820591655119\\
5.85858585858586	-0.165676483288665\\
6.06060606060606	-0.171767350700801\\
6.26262626262626	-0.177205146262581\\
6.46464646464646	-0.182081792382478\\
6.66666666666667	-0.186473366515618\\
6.86868686868687	-0.190443297824884\\
7.07070707070707	-0.194044832896823\\
7.27272727272727	-0.197322955426107\\
7.47474747474747	-0.200315894003747\\
7.67676767676768	-0.203056316264966\\
7.87878787878788	-0.205572282066079\\
8.08080808080808	-0.20788800994003\\
8.28282828282828	-0.210024497694718\\
8.48484848484848	-0.212000028201354\\
8.68686868686869	-0.21383058415569\\
8.88888888888889	-0.215530190173188\\
9.09090909090909	-0.217111196498876\\
9.29292929292929	-0.218584515517875\\
9.49494949494949	-0.219959819887226\\
9.6969696969697	-0.221245709289088\\
9.8989898989899	-0.222449851394062\\
10.1010101010101	-0.223579101522717\\
10.3030303030303	-0.224639604628926\\
10.5050505050505	-0.225636882546304\\
10.7070707070707	-0.226575908896912\\
10.9090909090909	-0.227461173628534\\
11.1111111111111	-0.228296738799484\\
11.3131313131313	-0.229086286949526\\
11.5151515151515	-0.229833163168488\\
11.7171717171717	-0.230540411789075\\
11.9191919191919	-0.231210808479221\\
12.1212121212121	-0.231846888385205\\
12.3232323232323	-0.232450970874249\\
12.5252525252525	-0.233025181340658\\
12.7272727272727	-0.233571470469258\\
12.9292929292929	-0.234091631291022\\
13.1313131313131	-0.234587314316949\\
13.3333333333333	-0.235060040994845\\
13.5353535353535	-0.235511215699196\\
13.7373737373737	-0.235942136434834\\
13.9393939393939	-0.236354004410493\\
14.1414141414141	-0.23674793261716\\
14.3434343434343	-0.237124953528291\\
14.5454545454545	-0.237486026023501\\
14.7474747474747	-0.237832041624418\\
14.9494949494949	-0.238163830119974\\
15.1515151515152	-0.238482164648799\\
15.3535353535354	-0.238787766298007\\
15.5555555555556	-0.239081308270487\\
15.7575757575758	-0.239363419666547\\
15.959595959596	-0.239634688920232\\
16.1616161616162	-0.239895666926128\\
16.3636363636364	-0.24014686988818\\
16.5656565656566	-0.240388781918476\\
16.7676767676768	-0.240621857410907\\
16.969696969697	-0.240846523211784\\
17.1717171717172	-0.241063180607044\\
17.3737373737374	-0.24127220714366\\
17.5757575757576	-0.241473958300902\\
17.7777777777778	-0.241668769025459\\
17.979797979798	-0.241856955143036\\
18.1818181818182	-0.242038814657597\\
18.3838383838384	-0.242214628948553\\
18.5858585858586	-0.242384663874774\\
18.7878787878788	-0.242549170793819\\
18.989898989899	-0.242708387503679\\
19.1919191919192	-0.242862539113727\\
19.3939393939394	-0.24301183885092\\
19.5959595959596	-0.243156488806683\\
19.7979797979798	-0.243296680629483\\
20	-0.243432596167487\\
};
\addlegendentry{$\alpha(\mu)$}

\end{axis}

\end{tikzpicture}%
    \end{minipage}
    \begin{minipage}{.48\linewidth}
      \input{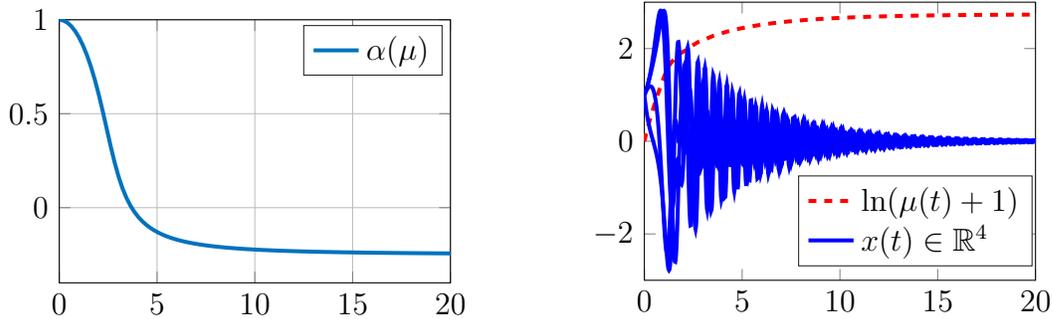}
    \end{minipage}
    \caption{Left: Spectral abscissa $\alpha_j$ as a function of
      $\mu$. Right: Adaptively stabilized
      system \eqref{eq:adaptively} with $x(0)=[1,1,1,1]^T$, $\mu(0)=0$.}\label{fig:specAbscHam}
  \end{figure}
\end{ex}

\subsection{Simultaneous stabilization by noise}
Stabilization of a dynamic system by noise processes is an interesting
phenomenon, which was analyzed in \cite{ArnoCrau83} (see also
e.g.\ \cite{Arno90, CaraRobi04}). As a particular situation, we consider the Stratonovich equation
\begin{align}\label{eq:StratoM}
  dx&=Ax\,dt+Mx\circ dw\;.
\end{align}
In this subsection we assume basic knowledge of stochastic calculus as e.g.\ in
\cite{Gard88a, KloePlat95}, but actually we only need the
spectral characterization of stability given in \eqref{eq:ItoLyap}.
Nevertheless we outline the background.
Informally, \eqref{eq:StratoM}  can be regarded as an ordinary differential
equation with noise perturbed coefficients, $\dot x(t)=(A+M\dot w(t))x(t)$.
Here $w(t)$ is a (stochastic) Wiener process, and the equation is understood as an
integral equation $x(t)=\int^t A x(\tau)\,d\tau+\int^t
Mx(\tau)\diamond dw(\tau)$ (the symbol $\diamond$ is explained below). The stochastic integral is approximated by
Riemann-Stieltjes type sums 
\begin{align*}
\int^t Mx(\tau)\diamond dw(\tau)=\lim\sum_{j=1}^N Mx(\tilde\tau_j)\big(w(\tau_j)-w(\tau_{j-1})\big)\;.
\end{align*}
Since  $w$ is not of bounded variation, the choice of
$\tilde\tau_j$ is essential. In the Stratonovich case (where we write $\diamond=\circ$), one sets
$\tilde\tau_j=(\tau_j+\tau_{j-1})/2$; in the It\^o-case (where
$\diamond$ is left out), one sets $\tilde\tau_j=\tau_j$.  While the
\emph{Stratonovich} interpretation is often more appropriate for \emph{modelling}
physical systems, \emph{analysis} and \emph{numerical solution} are easier for 
\emph{It\^o} equations. Therefore we will make use of transformations between
the solutions of the different types.
We call  \eqref{eq:StratoM}  \emph{asymptotically mean square} (or
\emph{2nd mean}) \emph{stable}, if for all solutions $x(t)$
the expected value of the squared norm $E(\|x(t)\|^2)$ converges to
zero as $t\to\infty$ (see e.g.\ \cite{KloePlat95, Damm04}).\\
For a given matrix $A\in\mathbb{R}^{n\times n}$ we want to construct
$M$ such that \eqref{eq:StratoM} is asymptotically mean square stable. 
It follows
from results in \cite{ArnoCrau83} that this is possible (with a
skew-symmetric $M$), if and only
if $\tr A<0$. Here we derive the following generalization.

\begin{theorem}\label{thm:ArnoCrau83plus}
Let $A_1,A_2\in\mathbb{R}^{n\times n}$ with $\tr A_1<0$ and $\tr
A_2<0$ be given. Then there exists a common skew-symmetric matrix $M$, such that
the systems
\begin{align}\label{eq:StratoMa}
  dx_1&=A_1x_1\,dt+Mx_1\circ dw_1\\
  dx_2&=A_2x_2\,dt+Mx_2\circ dw_2 \label{eq:StratoMb}
\end{align}
  are both asymptotically mean square stable. 
\end{theorem}

\bprf  
Let $\alpha_1=\tfrac{\tr A_1}n<0$ and $\alpha_2=\tfrac{\tr A_2}n<0$. 
By Proposition \ref{prop:fillmore_simultan} there exists an orthogonal
matrix $V$ such that $V^T(A_1-\alpha_1 I)V$ is hollow and
$V^T(A_2-\alpha_2 I)V$ is almost hollow. Transforming $x_j\mapsto V^Tx_j$ we can assume that  $A_1-\alpha_1 I$ is hollow and
$A_2-\alpha_2 I$ is almost hollow. 

For brevity we only elaborate on the case of  odd $n=2k+1$. The transfer to the
even case is then even easier (see Example \ref{ex:simstabnoise}). Let
$\omega=\left[\omega_1,\ldots,\omega_k\right]$ with
$0<\omega_1<\ldots<\omega_k$, and set
\begin{align*}
M(\omega)&=  \left[\begin{array}{cccc}
                    \begin{smallmatrix}
                      0
                    \end{smallmatrix}
&&&\\
    &\begin{smallmatrix}
     0& \omega_1\\-\omega_1&0
    \end{smallmatrix}
    &&\\
&&\ddots&\\
&&&\begin{smallmatrix}
     0& \omega_k\\-\omega_k&0
    \end{smallmatrix}
  \end{array}\right]\in\mathbb{R}^{n\times n}\;.
\end{align*}
We claim, that for $M=\mu M(\omega)$ with sufficiently large $\mu>0$
both \eqref{eq:StratoMa} and \eqref{eq:StratoMb} are asymptotically
mean square stable.\\
Note, that all eigenvalues of $M(\omega)$ are simple. An orthonormal set
of eigenvectors is given by $u_1=e_1$ and
$u_j=\frac1{\sqrt2}(e_j+ie_{j+1})$,
$u_{j+1}=\frac1{\sqrt2}(e_j-ie_{j+1})$ for even $j$.  
Hence with $U=[u_1,\ldots,u_n]$, we have
\begin{align}\label{eq:specdecM}
  U^*M(\omega)U=\diag(0,i\omega_1,-i\omega_1,\ldots,i\omega_k,-i\omega_k)=:\diag(i\tilde
  \omega_1,\ldots,i\tilde
  \omega_n)\;.
\end{align}

We rewrite the Stratonovich equations as the equivalent It\^o
equations (e.g.\ \cite{Gard88a})
\begin{align}\label{eq:ItoM}
  dx_j&=\left(A_j+\frac12M^2\right)x_j\,dt+Mx_j\,dw_j\;.
\end{align}
It is well-known (e.g.\ \cite{Damm04}), that \eqref{eq:ItoM} is asymptotically mean square stable, if and only if
\begin{align}\label{eq:ItoLyap}
  \sigma(\mathcal{L}_{A_j+\tfrac12M^2}+\Pi_M)\subset\mathbb{C_-}\;.
\end{align}
Here $\mathcal{L}_{N}:X\mapsto NX+XN^T$ for arbitrary
$N\in\mathbb{R}^{n\times n}$, and $\Pi_{M}:X\mapsto MXM^T$. We replace
$M$ by $\mu M(\omega)$. Then for large $\mu^2=1/\varepsilon$, we interpret
\begin{align}
\tfrac{1}{\mu^2}  \left(\mathcal{L}_{A_j+\tfrac12(\mu M(\omega))^2}+\Pi_{\mu M(\omega)}\right)
&=\left(\mathcal{L}_{M(\omega)^2/2}+\Pi_{M(\omega)}\right)+\varepsilon\mathcal{L}_{A_j}\label{eq:perturbedL}
\end{align}
as a perturbation of  $\mathcal{L}_{M(\omega)^2/2}+\Pi_{M(\omega)}$. It follows from \eqref{eq:specdecM} that
\begin{align*}
  (\mathcal{L}_{M(\omega)^2/2}+\Pi_M(\omega))(u_ku_\ell^*)&=\tfrac12\left(M(\omega)^2
                                            u_ku_\ell^*+u_ku_\ell^*M(\omega)^2\right)+M(\omega)
                                            u_ku_\ell^*M(\omega)\\
&=-\tfrac12\left(\tilde\omega_k^2+\tilde\omega_\ell^2-2\tilde\omega_k\tilde\omega_\ell\right) u_ku_\ell^* =-\tfrac12\left(\tilde\omega_k-\tilde\omega_\ell\right)^2 u_ku_\ell^*
\end{align*}
with $\tilde\omega_k-\tilde\omega_\ell=0$, if and only if $k=\ell$.
Thus, $\mathcal{L}_{M(\omega)^2/2}+\Pi_M(\omega)$ has an $n$-fold
eigenvalue~$0$ while all other eigenvalues are strictly negative.
We only have to consider the perturbation of the eigenvalue $0$. For
small $\varepsilon$, the perturbed mapping \eqref{eq:perturbedL} has an
$n$-dimensional invariant subspace with a basis, which depends
smoothly on $\varepsilon$ and coincides with $u_1u_1^*,\ldots,u_nu_n^*$
for $\varepsilon=0$, see \cite{StewSun90}. The restriction of   \eqref{eq:perturbedL} to this
subspace has the matrix representation $B_j=(b_{k\ell}^{(j)})$ with
\begin{align*}
  b_{k\ell}^{(j)}&=\tr\left(\mathcal{L}_{A_j}(u_\ell u_\ell^*\right)u_ku_k^*)=u_k^*\left(A_j
             u_\ell u_\ell^*+u_\ell u_\ell^*A_j^T\right)u_k\\&=\left\{
             \begin{array}{ll}
               0&\ell\neq k\\
              u_k^*(A_j+A_j^T)u_k=2\alpha_j&\ell= k
             \end{array}
\right. \;,
\end{align*}
since both $A_j-\alpha_j I$ are almost hollow. 
Hence
$B_j=2\alpha_j I$ has all eigenvalues in $\mathbb{C}_-$ and so has the
matrix in \eqref{eq:perturbedL} for sufficiently small $\varepsilon$.
This proves that for $M=\mu M(\omega)$ with sufficiently large $\mu$,
both 
  \eqref{eq:StratoMa} and \eqref{eq:StratoMb} are asymptotically mean square stable.
\eprf

\begin{ex}\rm\label{ex:simstabnoise}
  For an illustration with even $n$, we choose the simple but
  arbitrary matrix pair
  \begin{align*}
    (A_1,A_2)&=\left(\left[\begin{smallmatrix}
-1&1&1&1&1&1\\
  1&0&1&1&1&1\\
  0&1&0&1&1&1\\
  0&0&1&0&1&1\\
  0&0&0&1&0&1\\
  0&0&0&0&1&0 \end{smallmatrix}\right],
             \left[\begin{smallmatrix}1 & -1 & 0 & 0 & 0 & 0\\ 1 & 1 & -1 & 0 &
                 0 & 0\\ 1 & 0 & 1 & -1 & 0 & 0\\ 1 & 0 & 0 & 1 & -1 & 0\\ 1 & 0 & 0 & 0 & 1 & -1\\ 1 & 0 & 0 & 0 & 0 & -6 \end{smallmatrix}\right]\right)\;.
  \end{align*}
The  orthogonal matrix 
  \begin{align*}
U=    \left[\begin{smallmatrix} 
\phantom-0.1919&\phantom-0.1709&-0.1182&\phantom-0.4410&\phantom-0.3961&\phantom-0.7541\\
 -0.8960&-0.1266&\phantom-0.1726&-0.0363&-0.1203&\phantom-0.3682\\
  \phantom-0.0159&-0.6560&-0.1059&-0.3989&\phantom-0.6311&\phantom-0.0298\\
  \phantom-0.0144&\phantom-0.0086&-0.8175&-0.3556&-0.3660&\phantom-0.2664\\
  \phantom-0.0138&\phantom-0.6274&\phantom-0.2379&-0.6786&\phantom-0.2555&\phantom-0.1542\\
 -0.3996&\phantom-0.3616&-0.4692&\phantom-0.2411&\phantom-0.4808&-0.4473
\end{smallmatrix}\right]
  \end{align*}
transforms $(A_1,A_2)$ to $(\tilde
A_1,\tilde A_2)$ with $(\tilde A_1+\frac1nI,\tilde A_2+\frac1nI)$
being almost hollow, where
\begin{align*}
\tilde A_1&=\left[\begin{smallmatrix}
-0.1667&-0.6778&\phantom-0.8432&\phantom-0.5969&-1.2359&-0.8144\\
 \phantom-0.3655&-0.1667&-0.1359&\phantom-0.0294&-0.0818&-0.4453\\
  \phantom-0.4809&-0.4877&-0.1667&\phantom-1.1305&-1.0531&\phantom-0.2396\\
  \phantom-0.2712&-0.5650&\phantom-1.0652&-0.1667&-0.4391&-0.0790\\
 -1.3083&\phantom-0.7799&-0.8330&-1.1435&-0.1667&\phantom-0.0971\\
 -1.3506&\phantom-0.1132&-1.5411&-1.4969&\phantom-1.1554&-0.1667
\end{smallmatrix}\right]\;,\\
\tilde A_2&=\left[\begin{smallmatrix} 
-0.1667&\phantom-0.2200&-1.2765&-0.2157&\phantom-1.4333&-2.2393\\
 \phantom-1.4680&-0.1667&\phantom-0.8754&-0.9385&-1.5753&\phantom-1.6896\\
 -1.5017&\phantom-1.5458&-0.1667&-0.2265&\phantom-1.1226&-1.9108\\
\phantom-0.5741&-0.3164&\phantom-0.2973&-0.1667&-0.9509&-0.4748\\
\phantom-1.9688&-0.9634&\phantom-2.1166&-0.5422&\phantom-0.0562&\phantom-2.0303\\
 -0.4528&\phantom-1.3096&-1.5708&\phantom-1.2474&\phantom-1.3797&-0.3895
\end{smallmatrix}\right]\;. 
  \end{align*}
For
$
M\left(\left[
  \begin{smallmatrix}
    1\\2\\3
  \end{smallmatrix}
\right]\right)=  \left[\begin{array}{ccc}
    \begin{smallmatrix}
     0& 1\\-1&0
    \end{smallmatrix}
    &&\\
& \begin{smallmatrix}
     0& 2\\-2&0
    \end{smallmatrix}&\\
&&\begin{smallmatrix}
     0&3\\-3&0
    \end{smallmatrix}
  \end{array}\right]
$ 
we obtain the stabilizing skew-symmetric matrix
\begin{align*}
M=  UM\left(\left[
  \begin{smallmatrix}
    1\\2\\3
  \end{smallmatrix}
\right]\right)U^T&=\left[
                   \begin{smallmatrix}
\phantom-0.0000&\phantom-0.6949&-1.3331&\phantom-1.9489&-0.3262&-1.1247\\
 -0.6949&-0.0000&-0.2634&\phantom-0.1201&-1.1153&-0.6950\\
 \phantom-1.3331&\phantom-0.2634&\phantom-0.0000&-0.0300&\phantom-0.6217&-1.5717\\
 -1.9489&-0.1201&\phantom-0.0300&\phantom-0.0000&\phantom-0.9140&-0.6124\\
 \phantom-0.3262&\phantom-1.1153&-0.6217&-0.9140&\phantom-0.0000&-0.8317\\
\phantom-1.1247&\phantom-0.6950&\phantom-1.5717&\phantom-0.6124&\phantom-0.8317&-0.0000
                   \end{smallmatrix}
\right]\;.
\end{align*}
In Fig.\ \ref{fig:SpecAbsc}, we have plotted the spectral abscissae
\begin{align*}
\alpha_j(\mu)=\max\Real\sigma\left(\mathcal{L}_{A_j+\tfrac12(\mu M)^2}+\Pi_{\mu M}\right)
\end{align*}
for $j=1,2$ depending on $\mu$. Roughly for $\mu\ge 7$ both are
negative. We
chose $\mu=5$ and $\mu=20$ for simulations,  where
$\alpha_1(5)\approx-0.03<0$, $\alpha_2(5)=0.25>0$,
$\alpha_1(20)\approx-0.32<0$, $\alpha_2(20)=-0.29<0$. For both cases, 
Fig.\  \ref{fig:Samples} shows five sample paths of
$\|x_j\|$, $j=1,2$, with random initial conditions $x_0$ satisfying $\|x_0\|=1$. The solutions were computed by the
Euler-Maruyama scheme (e.g.\ \cite{KloePlat95}) with step size $1e-5$ applied to the It\^o
formulation \eqref{eq:ItoM} of the Stratonovich equation. The plots
exhibit the expected stability behaviour.
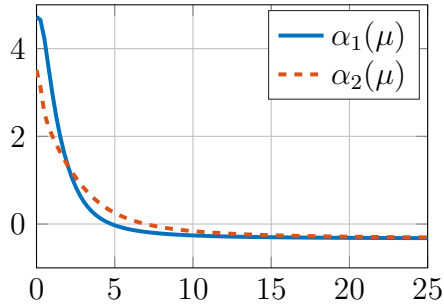
\begin{figure}[h]\centering
    \begin{minipage}{.5\linewidth}
%
%
\definecolor{mycolor1}{rgb}{0.00000,0.44700,0.74100}%
\definecolor{mycolor2}{rgb}{0.85000,0.32500,0.09800}%
\begin{tikzpicture}

\begin{axis}[%
width=5.2cm,
height=3.5cm,
at={(0cm,0cm)},
scale only axis,
xmin=0,
xmax=25,
ymin=-1,
ymax=5,
axis background/.style={fill=white},
xmajorgrids,
ymajorgrids,
legend style={legend cell align=left, align=left, draw=white!15!black}
]
\addplot [color=mycolor1, line width=1.5pt]
  table[row sep=crcr]{%
0	4.72467966571489\\
0.252525252525253	4.65492728530957\\
0.505050505050505	4.22484907393679\\
0.757575757575758	3.59952271507859\\
1.01010101010101	2.98431164617355\\
1.26262626262626	2.44686523116095\\
1.51515151515152	1.99516512601152\\
1.76767676767677	1.61970803255703\\
2.02020202020202	1.30760922341663\\
2.27272727272727	1.04775123039255\\
2.52525252525253	0.831800659465005\\
2.77777777777778	0.653273597293294\\
3.03030303030303	0.506502210265725\\
3.28282828282828	0.386232204581244\\
3.53535353535354	0.287676041360422\\
3.78787878787879	0.206665486901172\\
4.04040404040404	0.139723172547248\\
4.29292929292929	0.0840354408280754\\
4.54545454545455	0.0373675794469188\\
4.7979797979798	-0.00203886565446499\\
5.05050505050505	-0.0355636848080584\\
5.3030303030303	-0.0642911063536476\\
5.55555555555556	-0.0890766291033244\\
5.80808080808081	-0.110599152617318\\
6.06060606060606	-0.129401007620868\\
6.31313131313131	-0.145918486468099\\
6.56565656565657	-0.160505095901076\\
6.81818181818182	-0.173449301953297\\
7.07070707070707	-0.184988126667297\\
7.32323232323232	-0.195317621890879\\
7.57575757575758	-0.204600986391929\\
7.82828282828283	-0.212974897166265\\
8.08080808080808	-0.220554480474344\\
8.33333333333333	-0.227437240679731\\
8.58585858585859	-0.233706185638495\\
8.83838383838384	-0.239432328778026\\
9.09090909090909	-0.244676704517524\\
9.34343434343434	-0.249492001360821\\
9.5959595959596	-0.253923892748035\\
9.84848484848485	-0.25801212756325\\
10.1010101010101	-0.261791428395631\\
10.3535353535354	-0.265292235169988\\
10.6060606060606	-0.268541323726573\\
10.8585858585859	-0.271562322747118\\
11.1111111111111	-0.274376147632674\\
11.3636363636364	-0.277001366209257\\
11.6161616161616	-0.279454508210318\\
11.8686868686869	-0.28175032819426\\
12.1212121212121	-0.283902029721566\\
12.3737373737374	-0.285921457179139\\
12.6262626262626	-0.287819260477083\\
12.8787878787879	-0.28960503691274\\
13.1313131313131	-0.291287453745933\\
13.3838383838384	-0.292874354431958\\
13.6363636363636	-0.294372850949648\\
13.8888888888889	-0.295789404270448\\
14.1414141414141	-0.297129894677393\\
14.3939393939394	-0.298399683367888\\
14.6464646464646	-0.299603666564776\\
14.8989898989899	-0.30074632314611\\
15.1515151515152	-0.301831756674952\\
15.4040404040404	-0.30286373257104\\
15.6565656565657	-0.303845711045532\\
15.9090909090909	-0.304780876358344\\
16.1616161616162	-0.305672162842294\\
16.4141414141414	-0.306522278112508\\
16.6666666666667	-0.30733372379255\\
16.9191919191919	-0.308108814063071\\
17.1717171717172	-0.308849692290736\\
17.4242424242424	-0.309558345951926\\
17.6767676767677	-0.310236620062619\\
17.9292929292929	-0.310886229275764\\
18.1818181818182	-0.311508768790843\\
18.4343434343434	-0.312105724216565\\
18.6868686868687	-0.312678480490557\\
18.9393939393939	-0.313228329972689\\
19.1919191919192	-0.313756479773163\\
19.4444444444444	-0.314264058425827\\
19.6969696969697	-0.314752121952773\\
19.9494949494949	-0.31522165938863\\
20.2020202020202	-0.315673597825095\\
20.4545454545455	-0.316108807009516\\
20.7070707070707	-0.316528103554633\\
20.959595959596	-0.316932254793859\\
21.2121212121212	-0.317321982307347\\
21.4646464646465	-0.317697965160205\\
21.7171717171717	-0.318060842881051\\
21.969696969697	-0.318411218188684\\
22.2222222222222	-0.318749659512536\\
22.4747474747475	-0.319076703310012\\
22.7272727272727	-0.319392856197592\\
22.979797979798	-0.319698596924676\\
23.2323232323232	-0.319994378198102\\
23.4848484848485	-0.320280628362557\\
23.7373737373737	-0.320557752957353\\
23.989898989899	-0.320826136161123\\
24.2424242424242	-0.321086142125739\\
24.4949494949495	-0.321338116223625\\
24.7474747474747	-0.321582386193207\\
25	-0.321819263211229\\
};
\addlegendentry{$\alpha_1(\mu)$}

\addplot [color=mycolor2, dashed, line width=1.5pt]
  table[row sep=crcr]{%
0	3.51734350130203\\
0.252525252525253	3.12793128876446\\
0.505050505050505	2.52468011659331\\
0.757575757575758	2.22682538295092\\
1.01010101010101	2.00688627879764\\
1.26262626262626	1.8146136198829\\
1.51515151515152	1.63827592042292\\
1.76767676767677	1.47417800681147\\
2.02020202020202	1.32162913559323\\
2.27272727272727	1.18041778901823\\
2.52525252525253	1.05013939384095\\
2.77777777777778	0.930326049131636\\
3.03030303030303	0.820557195600696\\
3.28282828282828	0.720428333230164\\
3.53535353535354	0.629487723756543\\
3.78787878787879	0.547205249031213\\
4.04040404040404	0.47297771631975\\
4.29292929292929	0.40615442655089\\
4.54545454545455	0.346067869592717\\
4.7979797979798	0.292060779137508\\
5.05050505050505	0.243506149018598\\
5.3030303030303	0.199819976867132\\
5.55555555555556	0.160467923125219\\
5.80808080808081	0.12496746129116\\
6.06060606060606	0.0928869834502236\\
6.31313131313131	0.0638430317641027\\
6.56565656565657	0.037496511655387\\
6.81818181818182	0.0135484720313104\\
7.07070707070707	-0.00826417032616898\\
7.32323232323232	-0.0281727348911635\\
7.57575757575758	-0.0463805852879312\\
7.82828282828283	-0.0630665834606677\\
8.08080808080808	-0.0783881748592778\\
8.33333333333333	-0.0924841145499761\\
8.58585858585859	-0.10547685802864\\
8.83838383838384	-0.117474647147177\\
9.09090909090909	-0.128573323734533\\
9.34343434343434	-0.138857902991407\\
9.5959595959596	-0.148403936788621\\
9.84848484848485	-0.157278694325217\\
10.1010101010101	-0.165542184664672\\
10.3535353535354	-0.173248042759179\\
10.6060606060606	-0.180444297814016\\
10.8585858585859	-0.187174040332103\\
11.1111111111111	-0.193476001940368\\
11.3636363636364	-0.199385060114102\\
11.6161616161616	-0.204932678207312\\
11.8686868686869	-0.210147289696646\\
12.1212121212121	-0.215054634279711\\
12.3737373737374	-0.219678052362531\\
12.6262626262626	-0.224038743538096\\
12.8787878787879	-0.228155993860098\\
13.1313131313131	-0.232047376025398\\
13.3838383838384	-0.23572892601257\\
13.6363636363636	-0.239215299206443\\
13.8888888888889	-0.24251990863914\\
14.1414141414141	-0.245655047602026\\
14.3939393939394	-0.24863199857494\\
14.6464646464646	-0.25146113016714\\
14.8989898989899	-0.254151983522182\\
15.1515151515152	-0.256713349467094\\
15.4040404040404	-0.259153337488856\\
15.6565656565657	-0.261479437519346\\
15.9090909090909	-0.263698575343358\\
16.1616161616162	-0.265817162375105\\
16.4141414141414	-0.26784114043233\\
16.6666666666667	-0.269776022076746\\
16.9191919191919	-0.27162692699691\\
17.1717171717172	-0.273398614883022\\
17.4242424242424	-0.275095515155485\\
17.6767676767677	-0.276721753897378\\
17.9292929292929	-0.278281178275387\\
18.1818181818182	-0.279777378712383\\
18.4343434343434	-0.281213709049861\\
18.6868686868687	-0.282593304893732\\
18.9393939393939	-0.283919100339495\\
19.1919191919192	-0.285193843222334\\
19.4444444444444	-0.286420109052772\\
19.6969696969697	-0.287600313752478\\
19.9494949494949	-0.288736725315881\\
20.2020202020202	-0.289831474493184\\
20.4545454545455	-0.290886564586332\\
20.7070707070707	-0.291903880448377\\
20.959595959596	-0.292885196750322\\
21.2121212121212	-0.293832185583604\\
21.4646464646465	-0.294746423469536\\
21.7171717171717	-0.295629397804577\\
21.969696969697	-0.296482512818537\\
22.2222222222222	-0.297307095073822\\
22.4747474747475	-0.298104398536334\\
22.7272727272727	-0.298875609281946\\
22.979797979798	-0.299621849834542\\
23.2323232323232	-0.300344183203135\\
23.4848484848485	-0.301043616610792\\
23.7373737373737	-0.301721104959477\\
23.989898989899	-0.302377554047843\\
24.2424242424242	-0.303013823560923\\
24.4949494949495	-0.303630729848784\\
24.7474747474747	-0.304229048516227\\
25	-0.304809516827205\\
};
\addlegendentry{$\alpha_2(\mu)$}

\end{axis}

\begin{axis}[%
width=6.135cm,
height=4.294cm,
at={(-0.798cm,-0.472cm)},
scale only axis,
xmin=0,
xmax=1,
ymin=0,
ymax=1,
axis line style={draw=none},
ticks=none,
legend style={legend cell align=left, align=left, draw=white!15!black}
]
\end{axis}
\end{tikzpicture}%
    \end{minipage}
\caption{Spectral abscissa $\alpha_j$ as a function of
  $\mu$.}\label{fig:SpecAbsc}
  \end{figure}
\begin{figure}[h]
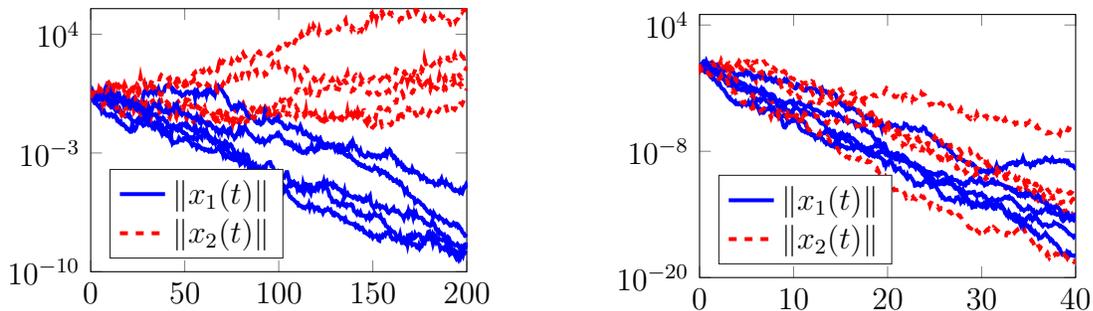
\centering
    \begin{minipage}{.5\linewidth}
      \input{./SimStrato5a.tex}
    \end{minipage}\hfill
 \begin{minipage}{.5\linewidth}
      \input{./SimStrato.tex}
    \end{minipage}
\caption{Sample paths of $\|x_j(t)\|$ for $\mu=5$ (left) and
  $\mu=20$ (right)}\label{fig:Samples}
  \end{figure}
\end{ex}
\begin{remark}\label{rem:common22}
  There even exists a common skew-symmetric matrix $M$ so that $m$ equations
\begin{align}\label{eq:StratoMj}
  dx_j&=A_jx_j\,dt+Mx_j\circ dw_j \text{ with }\tr A_j<0  \quad j=1,\ldots,m 
\end{align}
are simultaneously stabilized, if a common  orthogonal matrix $U$ can be
found, so that for all $j$ 
\begin{align*}
\diag\left(U^T\Big(A_j-\tfrac{\tr  A_j}nI\Big)U\right)&=[d^{(j)}_1,-d^{(j)}_1,\ldots,d^{(j)}_k,-d^{(j)}_k,0],
                       \text{ if } n=2k+1, \text{ or }\\
\diag\left(U^T\Big(A_j-\tfrac{\tr
  A_j}nI\Big)U\right)&=[d^{(j)}_1,-d^{(j)}_1,\ldots,d^{(j)}_k,-d^{(j)}_k],
                       \text{ if } n=2k.
\end{align*}
The proof of Theorem \ref{thm:ArnoCrau83plus} applies literally in
this case.\\
If the matrix $U$ can be chosen symplectic, then $M$ can be chosen
Hamiltonian, as a combination with the proof of Theorem
\ref{thm:HamSkew} shows. 
\end{remark}
\section{Conclusion and outlook}
As our main theoretic contribution we see Theorem \ref{thm:SymplOrth},
which states that every real matrix is symplectic-orthogonally similar to a
matrix with constant diagonal, (w.l.o.g.\ a hollow matrix, if the
trace is subtracted). The proof requires a result on the simultaneous
transformation of two matrices which is closely related to properties
of the joint numerical range. For our applications it turns out that
the hollow form can be weakened to a $2\times 2$-block hollow form,
where only $a_{ii}+a_{i+1,i+1}=0$ for $i=1,3,\ldots$ (see Remark \ref{rem:common22}).  
This gives rise to further connections and questions, which were not
discussed here. For instance, a simultaneous transformation to a
$2\times 2$-block hollow form is related to the real $2$-nd numerical range (cf.\
\cite{FillWill71, LiPoon00}). General conditions on the convexity of
the real $2$-nd numerical range (like e.g.\ in \cite{GutkJonc04}) do
not seem to be available. Therefore it is unclear, whether more than two
zero-trace  matrices can always be transformed to  $2\times 2$-block hollow form.\\
Numerically, also the following variant of Proposition
\ref{prop:fillmore_simultan} seems to hold, but we were not able to
prove it. We state it as a conjecture.
\begin{conj}
  Consider $A,B\in\mathbb{R}^{n\times n}$ with
  $\tr A=\tr B=0$.
  There exists an orthogonal matrix $V\in\mathbb{R}^{n\times
     n}$ such that $V^TAV$ is hollow and $VBV^T$ is almost hollow.
   Note that here $A\mapsto V^TAV$, but $B\mapsto VBV^T$; the
   transformation applied to $A$ is the inverse (also adjoint) of the
   one applied to $B$ (unlike in Proposition
\ref{prop:fillmore_simultan}). 
\end{conj}

\appendix

\section{Direct symplectic orthogonal transformation of a symmetric
  $4\times 4$ matrix}\label{app.1}
In this appendix we develop a much more efficient alternative to the
3rd step in the proof of Theorem~\ref{thm:SymplOrth}. Using an
adapted notation, we now consider the symmetric $4\times 4$-matrix
\begin{align}
A = \begin{bmatrix}
a & b & c & d\\
b & e & f & g\\
c & f & h & i\\
d & g & i & j
\end{bmatrix} \quad \text{with}\quad a+e+h+j = 0. \label{eq:Aabcd}
\end{align}
\begin{claim}
  There exists a symplectic orthogonal transformation
  $S\in\mathbb{R}^{4\times 4}$, such that
  $S^TAS$ for $A$ in \eqref{eq:Aabcd} is hollow.
\end{claim}
\bprf
  We will obtain $S$ as the product of two symplectic orthogonal transformations
  $S=S_1S_2$ and consider $S^TAS = S_2^T\left(S_1^TAS_1\right)S_2$.
  Assume that for $\tilde{A}= S_1^TAS_1$ we have
  $\tilde{a}_{11}=-\tilde{a}_{33}~~ (\tilde{a} = -\tilde{h})$ and
  $\tilde{a}_{22}=-\tilde{a}_{44}~~ (\tilde{e}=-\tilde{j}).$ Consider
  the symplectic orthogonal matrix
  \[
  S_2 = \begin{bmatrix}
    q_0 & 0 & q_2 & 0\\
    0 &q_1 & 0 & q_3\\
    -q_2 & 0 & q_0 & 0\\
    0 & -q_3 & 0 & q_1
  \end{bmatrix} \quad \text{with}\quad q_0^2+q_2^2= 1, q_1^2+q_3^2 = 1
  \]
  and its effect on the diagonal elements of $\tilde{A}.$ (Note that
  $S_2=G_1(q_0,q_2)G_2(q_1,q_3)$ with $G_1,G_2$ from
  \eqref{eq:G2nics}). That is, we consider the diagonal elements of
  $\hat{A}=S_2^T\tilde{A}S_2$
  \begin{align*}
    e_1^T\hat{A}e_1 &= q_0^2\tilde{a}-2q_0q_2\tilde{c}+q_2^2\tilde{h} = 0,\\
    e_2^T\hat{A}e_2 &=q_1^2\tilde{e}-2q_1q_3\tilde{g}+q_3^2\tilde{j} =0,\\
    e_3^T\hat{A}e_3 &=q_0^2\tilde{h}+2q_0q_2\tilde{c}+q_2^2\tilde{a} =0,\\
    e_4^T\hat{A}e_4 &=q_1^2\tilde{j}+2q_1q_3\tilde{g}+q_3^2\tilde{e} =0.
  \end{align*}
  The first and the third as well as the second and the fourth
  equation are identical as $\tilde{a} = -\tilde{h}$ and
  $\tilde{e}=-\tilde{j}.$ Thus it suffices to consider
  \begin{align*}
    e_1^T\hat{A}e_1&=q_0^2\tilde{a}-2q_0q_2\tilde{c}-q_2^2\tilde{a} = 0,\\
    e_2^T\hat{A}e_2 &=q_1^2\tilde{e}-2q_1q_3\tilde{g}-q_3^2\tilde{e} =0.
  \end{align*}
  As both equations have the same form, we only consider the first
  equation and divide by $q_0^2\tilde{a}$ (assuming that
  $\tilde{a}\neq0$)
  \[
  0= t^2 +2t\frac{\tilde{c}}{\tilde{a}}-1
  \]
  and
  \[
  t_{1,2} = -\frac{\tilde{c}}{\tilde{a}} \pm
  \sqrt{\left(\frac{\tilde{c}}{\tilde{a}}\right)^2+1}, \quad q_0 =
  \frac{1}{\sqrt{t^2+1}}, \quad q_2 = tq_0.
  \]
  In case $\tilde{a} = 0,$ the choice $q_0=1, q_2 = 0$ will be
  perfect. Thus, if $S_1$ can be chosen such that the
  diagonal elements of $\tilde{A}= S_1^TAS_1$ satisfy
  $\tilde{a}=-\tilde{h}$ and $\tilde{e}=-\tilde{j},$, then a symplectic
  orthogonal matrix $S_2$ can be constructed such that
  $\diag(S_2'\tilde{A}S_2) = 0.$

  Now let us consider
  \[
  S_1 = \begin{bmatrix}
    p_0 & -p_1 & -p_2& -p_3\\
    p_1 & p_0 & -p_3 & p_2\\
    p_2 & p_3 & p_0 & -p_1\\
    p_3 & -p_2 & p_1 & p_0
  \end{bmatrix}
  \]
  as in \eqref{eq:4x4sympl} and its effect on the diagonal elements of
  $\tilde{A}= S_1^TAS_1$
  \begin{align*}
    e_1^T\tilde{A}e_1
 &= p_0^2a+p_1^2e+p_2^2h+p_3^2j+
   2p_0(p_1b+p_2c+p_3d)+   
   2p_1(p_2f+p_3g)+2p_2p_3i, 
    \\
    e_2^T\tilde{A}e_2
 &=p_0^2e+p_1^2a+p_2^2j+p_3^2h+
   2p_0(-p_1b-p_2g+p_3f)  
   +2p_1(p_2d-p_3c)-2p_2p_3i
   ,\\
    e_3^T\tilde{A}e_3
 &= p_0^2h+p_1^2j+p_2^2a+p_3^2e
   +2p_0(p_1i-p_2c-p_3f)   
   +2p_1(-p_2d-p_3g)+2p_2p_3b 
   ,\\
    e_4^T\tilde{A}e_4
 &=p_0^2j+p_1^2h+p_2^2e+ p_3^2a
   +2p_0(-p_1i+p_2g-p_3d)-2p_1(p_2f-p_3c)  
   -2p_2p_3b.
  \end{align*}
  Now choose $p_0,p_1,p_2,p_3$ such that
  $e_1^T\tilde{A}e_1 =-e_3^T\tilde{A}e_3 $ and
  $e_2^T\tilde{A}e_2 =-e_4^T\tilde{A}e_4 $
  \begin{align*}
    &p_0^2a+p_1^2e+p_2^2h+p_3^2j+
      2p_0p_1b+2p_0p_2c+2p_0p_3d+
      2p_1p_2f+2p_1p_3g+2p_2p_3i\\
    &= -p_0^2h-p_1^2j-p_2^2a-p_3^2e
      -2p_0p_1i+2p_0p_2c+2p_0p_3f
      +2p_1p_2d+2p_1p_3g-2p_2p_3b,\\
    &p_0^2e+p_1^2a+p_2^2j+p_3^2h
      -2p_0p_1b-2p_0p_2g+2p_0p_3f
      +2p_1p_2d-2p_1p_3c-2p_2p_3i\\
    &=-p_0^2j-p_1^2h-p_2^2e- p_3^2a
      +2p_0p_1i-2 p_0p_2g+2p_0p_3d+2p_1p_2f
      -2p_1p_3c
      +2p_2p_3b,
  \end{align*}
  that is,
  \begin{align*}
    0&=(p_0^2+p_2^2)(a+h)+(p_1^2+p_3^2)(e+j)\\&\phantom=+
       2(p_0p_1+p_2p_3)(b+i)+2(p_0p_3-p_1p_2)(d-f),\\[1mm]
    0 &=(p_0^2+p_2^2)(e+j)+(p_1^2+p_3^2)(a+h)\\&\phantom=
       -2(p_0p_1+p_2p_3)(b+i)
       -2(p_0p_3-p_1p_2)(d-f)
  \end{align*}
  and
  \[p_0^2+p_1^2+p_2^2+p_3^2=1.
  \]
  Recall that $a+e+h+j=0$ holds. Thus, $a+h = -(e+j).$
  \begin{itemize}
  \item In case $a+h=e+j=0,$ we obtain the two equations
    \[
    0=(p_0p_1+p_2p_3)(b+i)+(p_0p_3-p_1p_2)(d-f), \quad
    p_0^2+p_1^2+p_3^2+p_3^2=1,
    \]
    which are satisfied for the choice $p_0=1, p_1=p_2=p_3=0,$ that
    is, $S_1 = I.$
  \item In case $a+h=-(e+j) \neq 0,$ we obtain the two equations
    \begin{align*}
      0&=(-p_0^2-p_2^2+p_1^2+p_3^2)(e+j)\\&+
         2(p_0p_1+p_2p_3)(b+i)+2(p_0p_3-p_1p_2)(d-f),\\[1mm]
      1&=p_0^2+p_1^2+p_2^2+p_3^2
    \end{align*}
    As $p_1^2+p_3^2 = 1-p_0^2-p_2^2$ this can be rewritten to
    \begin{align*}
      0&=(\tfrac12-p_0^2-p_2^2)(e+j)+
         (p_0p_1+p_2p_3)(b+i)+(p_0p_3-p_1p_2)(d-f),\\
      1&=p_0^2+p_1^2+p_2^2+p_3^2
    \end{align*}
    We need to distinguish some cases.
    \begin{itemize}
    \item In case $b+i = d-f=0$ we have $\tfrac12-p_0^2-p_2^2 =0,$
      that is $\tfrac12=p_0^2+p_2^2$ and $\tfrac12=p_1^2+p_3^2.$ One
      option is to choose $p_0=p_1=p_2=p_3 = \tfrac12.$ A different
      option is the choice $p_1=p_2 =0$ and
      $p_0 = p_3 = \frac{1}{\sqrt{2}},$ while a third option is given
      by $p_0=p_1 =0$ and $p_2 = p_3 = \frac{1}{\sqrt{2}}.$
    \item In case $b+i =0$ we have
      \begin{align*}
        0&=(\tfrac12-p_0^2-p_2^2)(e+j)+(p_0p_3-p_1p_2)(d-f),\\
        1&=p_0^2+p_1^2+p_2^2+p_3^2.
      \end{align*}
      The choice $p_1=p_3=0$ and $p_0 = p_2 = \frac{1}{\sqrt{2}}$
      yields the desired transformation.
    \item In case $d-f=0$ we have
      \begin{align*}
        0&=(\tfrac12-p_0^2-p_2^2)(e+j)+
           (p_0p_1+p_2p_3)(b+i),\\
        1&=p_0^2+p_1^2+p_2^2+p_3^2.
      \end{align*}
      The choice $p_1=p_3=0$ and $p_0 = p_2 = \frac{1}{\sqrt{2}}$
      yields the desired transformation.
    \item For all other cases, the simple choice $p_2 = p_3 = 0,$
      $S_1= \left[\begin{smallmatrix}
          p_0 & -p_1 & 0& 0\\
          p_1 & p_0 & 0 & 0\\
          0 & 0 & p_0 & -p_1\\
          0 & 0 & p_1 & p_0
        \end{smallmatrix}\right],$
      with $1=p_0^2+p_1^2$ and $0=(p_1^2-p_0^2)(e+j)+2p_0p_1(b+i)$
      gives the desired transformation. The second equation
      gives 
      with $t = \frac{p_0}{p_1}$
      \[t^2-2t\frac{b+i}{e+j}-1=0.\] Thus
      \[t_{1,2}=\frac{b+i}{e+j}\pm
      \sqrt{\left(\frac{b+i}{e+j}\right)^2+1}\] and
      \[ p_1 = \frac{1}{\sqrt{1+t^2}} \qquad p_0 = tp_1.\]
      Note that there are a number of other possible choices of
      $p_0, p_1,p_2,p_3,$ thus the symplectic orthogonal matrix $S_1$
      is not unique.
    \end{itemize}
  \end{itemize}
\eprf
\section*{Acknowledgments}
We thank Michael Karow for pointing out \cite{Bind85} to us,  Jean-Baptiste
Hiriart-Urruty for providing us with \cite{Pepi04}, and Richard Greiner
for sending us \cite{Yaku71}.


\end{document}